\documentclass[11pt,fleqn]{article}
\usepackage[utf8]{inputenc}
\usepackage{amsmath, amssymb, amsthm}
\usepackage{geometry}
\usepackage{hyperref}
\usepackage{microtype}
\usepackage{tabularx} 
\usepackage{xcolor}

\newcolumntype{C}{>{\centering\arraybackslash}X}

\newtheorem{theorem}{Theorem}
\newtheorem{lemma}{Lemma}
\newtheorem{definition}{Definition}
\newtheorem{proposition}{Proposition}
\newtheorem{corollary}{Corollary}
\theoremstyle{remark}
\newtheorem{remark}{Remark}

\title{\textbf{Unique Highest Trees of Random Mappings}}
\author{\textbf{Mikhail V. Berlinkov}
}
\date{\today}

\begin{document}

\maketitle

\begin{abstract}
We study the heights of the trees of a random mapping of $n$ elements. Using singularity analysis of exact generating functions, we prove that the mapping has a unique highest tree with probability $1 - \sqrt{\frac{\pi}{8}}\, n^{-1/2} + o(n^{-1/2})$. The property of having a unique highest tree plays a crucial role in the solution of the Road Coloring Problem \cite{trahtman2009}. More generally, we consider $c$-branches, the subtrees rooted at distance $c$ from the cycles. For all fixed $c \ge 0$ and $j \ge 1$ we show that the highest $c$-branch exceeds every other $c$-branch in height by at least $j$, except with probability asymptotic to $(2j-1) \sqrt{\frac{\pi}{8}}\, n^{-1/2}$. We also show that, for any fixed $\alpha > 1$, with probability $1 - \Theta(n^{-1/2})$ the part of the highest $c$-branch that lies above the second highest one (its crown) has more than $\alpha$ times as many vertices as roots. The last result is used in the author's proof that a random $2$-letter automaton with $n$ states is synchronizing with probability $1 - \Theta(1/n)$ \cite{berlinkov2016}.
\end{abstract}

\section{Introduction}
A random mapping is a function $f: \{1, \dots, n\} \to \{1, \dots, n\}$ chosen uniformly at random. Its functional graph is a digraph in which every vertex has out-degree $1$; it decomposes uniquely into weakly connected components, each consisting of a directed cycle of rooted trees.

Random mappings have been studied since the work of Harris \cite{harris1960} and Stepanov \cite{stepanov1969}. The heights of the trees involved are governed by the classical theory of the height of random trees: R\'{e}nyi and Szekeres \cite{renyi1967} showed that the height $H_n$ of a random labelled rooted tree with $n$ vertices has mean asymptotic to $\sqrt{2\pi n}$ and that $H_n/\sqrt{2n}$ has a limit distribution of theta type, and Flajolet and Odlyzko \cite{flajolet1982} extended the analysis of the average height to all simple varieties of trees by singularity analysis of the generating functions of trees of bounded height. The continuous limits of these structures are described by Aldous's continuum random tree \cite{aldous1991} and, for random mappings, by the Brownian bridge asymptotics of Aldous and Pitman \cite{aldous1994}; see also \cite{mutafchiev1988} and the monograph \cite{drmota2009}.

These results describe the height of the highest tree of a random mapping, which is of order $\sqrt{n}$. In this paper we are interested in a finer question: how the highest tree, or more generally the highest subtree rooted at a given distance from the cycles, compares with its competitors.

\paragraph{Results.} Our first result (Theorem~\ref{th_main}) is that a random mapping of $n$ elements has a unique highest tree with probability $1 - \sqrt{\pi/8}\; n^{-1/2} + o(n^{-1/2})$. For the application described below we need a more general and more robust statement. For an integer $c \ge 0$, a \emph{$c$-branch} is a subtree rooted at a vertex at distance $c$ from the cycles, so that the $0$-branches are the trees of the mapping. We show (Theorem~\ref{th_gap}) that for all fixed $c \ge 0$ and $j \ge 1$ the highest $c$-branch is unique and exceeds every other $c$-branch in height by at least $j$, except with probability asymptotic to $(2j-1) \sqrt{\pi/8}\; n^{-1/2}$; in particular, the constant $\sqrt{\pi/8}$ of the first result does not depend on $c$. Finally, we consider the \emph{crown} of the highest $c$-branch, that is, the part of it that lies above the height of the second highest $c$-branch, and we show (Theorem~\ref{th_high_tree}) that for every fixed $\alpha > 1$ the crown has more than $\alpha$ times as many vertices as roots with probability $1 - \Theta(n^{-1/2})$.

\paragraph{Method.} The proofs use the symbolic method and singularity analysis \cite{flajolet2009}. The events in question have exact generating functions, which are sums over the height $h$ of rational expressions in the generating functions $t_h(z)$ of trees of height at most $h$ (Sections~\ref{sec:unique_tree} and~\ref{sec:branches}). The analytic core of the paper is Lemma~\ref{lem:uniform}, which provides estimates for $t_h(z)$ that are uniform in $h$ and in a $\Delta$-domain at the singularity $z = e^{-1}$, and which shows that the sums have a logarithmic singularity with an explicit coefficient. The estimates are derived from the recurrence satisfied by $t - t_h$, as in \cite{flajolet1982}. For the crown we bound the generating function coefficientwise by an exponentially tilted one, which leads to a second solution of the same recurrence.

\paragraph{Motivation.} Our primary motivation comes from the theory of synchronizing automata. An automaton with $n$ states in our terminology can be defined as a finite set $A$ of mappings (called letters) on $n$ states, and one is called synchronizing if some composition of the letters from $A$ forms a constant mapping. For an introduction to the theory of synchronizing automata, we refer the reader to the survey \cite{volkov2008}.

The bridge between these two areas was established by Trahtman \cite{trahtman2009} in his breakthrough proof of the Road Coloring Problem: ``Given a strongly-connected digraph with constant out-degree $k>1$ and a greatest common divisor (gcd) of cycle lengths equal to one, is it always possible to find a coloring of all outgoing arrows for each vertex of the graph in $k$ colors which would result in a synchronizing automaton?'' The problem had been open for more than thirty years despite numerous attempts to solve it (see e.g., \cite{kari2002}, \cite{obrien1981}, \cite{kari2003}). The proof demonstrates that any valid graph can be colored such that the functional graph of one of the letters contains a unique highest tree, and shows that this configuration is sufficient to reduce the problem to a smaller size.

Another challenging but younger problem is about the probability that a random automaton with $n$ states is synchronizable. In \cite{CamConj} Cameron conjectured that the asymptotic probability for the most interesting case of 2-letter automata is $1-\Theta(\frac{1}{n})$ and some partial results were established in \cite{nicaud2016}, \cite{zaks2010}, \cite{zaks2013}. It turns out that a property similar to having a unique highest tree is useful for attacking this problem, and an asymptotic estimate confirming Cameron's conjecture was established in the author's work \cite{berlinkov2016}, conditional on the results proved here. To be able to leverage independence of the letters in~\cite{berlinkov2016}, we need a \emph{stable pair} of states to be completely defined by one letter with this extreme property. One state of this pair is defined as a first-level vertex in the (unique) highest tree. However, with this definition this state can be any of the first-level vertices in this tree and hence is not completely defined by the graph of this letter. This is the reason for considering $c$-branches instead of trees; the application uses the case $c = 1$ of Theorem~\ref{th_high_tree}.

\paragraph{Earlier versions.} Earlier versions of this manuscript (arXiv:1504.04532, versions 1--7) approached these results through the theory of Galton--Watson branching processes. That argument relied on details of proofs in the literature which were difficult to verify, and it gave the constant $\frac{2\pi}{3} - \frac{827}{288\pi} \approx 1.18$ in place of $\sqrt{\pi/8} \approx 0.63$ in Theorem~\ref{th_main}. That value is incorrect, as the exact values in Section~\ref{sec:empitrical} also show. The present version is based on a different, self-contained method. The application in \cite{berlinkov2016} only uses the order $\Theta(n^{-1/2})$ and is not affected by the change of the constant.

\paragraph{Organization.} Section~\ref{sec:prelim} recalls the tree function and the transfer theorems. Section~\ref{sec:unique_tree} contains the generating function of mappings without a unique highest tree, the uniform estimates (Lemma~\ref{lem:uniform}, stated with an outline of its proof in Section~\ref{sec:uniform}) and the proof of Theorem~\ref{th_main}. Section~\ref{sec:branches} treats $c$-branches: exact generating functions, the gap between the two highest $c$-branches (Theorem~\ref{th_gap}), and the crown (Theorem~\ref{th_high_tree}). Section~\ref{sec:empitrical} reports exact numerical values and simulations. The proof of Lemma~\ref{lem:uniform} is given in Appendix~\ref{app:uniform}, and the analytic part of the proof of Theorem~\ref{th_high_tree} in Appendix~\ref{app:crown}.

\section{Analytic Preliminaries}
\label{sec:prelim}

Let $\mathcal{T}$ be the class of standard Cayley trees. By the exponential formula for labeled sets and cycles, the Exponential Generating Function (EGF) for the class of all unconstrained mappings $\mathcal{M}$ is given by:
\begin{equation}
M(z) = \exp\left( \ln \left( \frac{1}{1-t(z)} \right) \right) = \frac{1}{1-t(z)}
\end{equation}

We rely on the following property of the tree function \cite[Proposition VI.3]{flajolet2009}:

\begin{theorem}[Singularity of Cayley Trees, Flajolet and Sedgewick 2009]\label{thm:tree}
The structural tree function $t(z)$ defined by the functional equation $t(z) = z \exp(t(z))$ has a dominant singularity at $z = e^{-1}$. In a $\Delta$-domain originating at $z = e^{-1}$, the function possesses the local expansion:
\begin{equation}
t(z) = 1 - y + \mathcal{O}(y^2)
\end{equation}
where $y = \sqrt{2(1 - ez)}$ is the local scaling variable.
\end{theorem}

This implies the mapping function scales as $M(z) \sim y^{-1} = (2(1-ez))^{-1/2}$. To extract asymptotic probabilities, we utilize the transfer theorems for algebraic and logarithmic singularities \cite{flajolet1990}:

\begin{theorem}[Transfer for Algebraic and Logarithmic Singularities, Flajolet and Odlyzko 1990]
\label{thm:transfer}
Let $A(z)$ be a generating function analytic in a $\Delta$-domain originating at its dominant singularity $z = \rho$.
\begin{itemize}
\item If $A(z) \sim \left(1 - \frac{z}{\rho}\right)^{-\alpha}$ as $z \to \rho$ for $\alpha \notin \{0, -1, -2, \dots\}$, then:
\begin{equation}
[z^n] A(z) \sim \rho^{-n} \frac{n^{\alpha - 1}}{\Gamma(\alpha)}
\end{equation}
\item If $A(z) = \lambda \ln \left( \frac{1}{1 - z/\rho} \right) + \kappa + o(1)$ as $z \to \rho$ in the $\Delta$-domain, for constants $\lambda$ and $\kappa$, then:
\begin{equation}
[z^n] A(z) = \lambda\, \frac{\rho^{-n}}{n} + o\left( \frac{\rho^{-n}}{n} \right)
\end{equation}
\item If $A(z) - \lambda \ln \left( \frac{1}{1 - z/\rho} \right)$ is bounded in the $\Delta$-domain for a constant $\lambda$, then:
\begin{equation}
[z^n] A(z) = \mathcal{O}\left( \frac{\rho^{-n}}{n} \right)
\end{equation}
\end{itemize}
\end{theorem}

Applying Theorem \ref{thm:transfer} to $M(z)$, the algebraic singularity maps to an exact asymptotic coefficient growth of $[z^n] M(z) \sim \frac{e^n}{\sqrt{2\pi n}}$.

\section{The Unique Highest Tree (0-Branch)}
\label{sec:unique_tree}

\begin{theorem}\label{th_main}
Let $\Sigma_n$ denote the set of all mappings on $n$ elements. The probability that a random mapping $g \in \Sigma_n$ possesses a unique highest tree is $1 - \sqrt{\frac{\pi}{8}}\, n^{-1/2} + o(n^{-1/2})$.
\end{theorem}

To track highest trees, we define $t_h(z)$ as the univariate EGF of trees bounded by height $h$. This satisfies the discrete recurrence $t_h(z) = z \exp(t_{h-1}(z))$ with the initial boundary condition $t_0(z) = z$ (representing a tree of maximum height 0 consisting solely of its root). The difference $\Delta_h(z) = t_h(z) - t_{h-1}(z)$ isolates the layer of trees achieving exactly height $h$.

To substitute height constraints into the cycle mapping operator, we define a bivariate generating function $t_h(a,z)$ for trees bounded by height $h$. Within this component, the parameter $a$ functions strictly as a formal marking parameter tracking whether the individual tree reaches exactly height $h$:
\begin{equation}
t_h(a,z) = t_{h-1}(z) + a(t_h(z) - t_{h-1}(z))
\end{equation}

\begin{lemma}
The bivariate EGF for functional mappings whose maximum tree height is bounded by $h$, marked by $a$ for the total multiplicity of trees achieving height $h$ across the mapping, is $M_h(a, z) = (1 - t_h(a,z))^{-1}$.
\end{lemma}

\begin{proof}
Substituting the bivariate component $t_h(a,z)$ into the set-of-cycles structural operator isolates the global maximum height, allowing the exponent of $a$ to track the total multiplicity of maximum-height trees.
\end{proof}

By extracting coefficients with respect to $a$, the configuration where exactly one tree reaches height $h$ is given by:
\begin{equation}
[a^1] M_h(a, z) = \left. \frac{\partial}{\partial a} \left[ \frac{1}{1 - t_h(a,z)} \right] \right|_{a=0} = \frac{\Delta_h(z)}{(1 - t_{h-1}(z))^2}
\end{equation}
The configuration where a structural tie occurs at maximum height $h$ is $M_{\mathrm{ties}, h}(z) = M_h(1, z) - M_h(0, z) - [a^1] M_h(a, z)$. Simplifying this yields:
\begin{equation}\label{eq:ties_discrete}
M_{\mathrm{ties}, h}(z) = \frac{(\Delta_h(z))^2}{(1 - t_{h-1}(z))^2 (1 - t_h(z))}
\end{equation}

Given the generating function for a tie at a specific height $h$, $M_{\mathrm{ties},h}(z)$, we define the global generating function for tie configurations across all possible heights as the infinite sum:
\begin{equation}
M_{\mathrm{ties}}(z) = \sum_{h=1}^{\infty} M_{\mathrm{ties},h}(z)
\end{equation}
The series $M_{\mathrm{ties}}(z)$ accounts for the ties at heights $h \ge 1$. The only other mappings without a unique highest tree are those in which all trees have height $0$ and there are at least two of them, that is, the $n!$ permutations of $n \ge 2$ elements. Since $|\Sigma_n| = n^n = n!\, [z^n] M(z)$, we have for $n \ge 2$
\begin{equation}\label{eq:tieprob}
\mathbb{P}(\text{$g$ has no unique highest tree}) = \frac{[z^n] M_{\mathrm{ties}}(z) + 1}{[z^n] M(z)} .
\end{equation}
The rest of this section determines the singular behaviour of $M_{\mathrm{ties}}(z)$ at $z = e^{-1}$. The estimates are established in a form that also covers the sums occurring in Section~\ref{sec:branches}.

\subsection{Uniform estimates for trees of bounded height}
\label{sec:uniform}

\paragraph{Notation.} We use the conventions $t_{-1}(z) := 0$ and $\Delta_0(z) := t_0(z) = z$, and we put
\[
\varepsilon_h(z) := 1 - t_h(z) \quad (h \ge -1), \qquad L(z) := \ln \frac{1}{1 - ez},
\]
with the principal branch of the logarithm. For $\eta > 0$ and $0 < \phi < \pi/2$ let
\[
\mathcal{D}(\eta, \phi) := \bigl\{ z : |z| < (1 + \eta) e^{-1},\ z \ne e^{-1},\ |\arg(z - e^{-1})| > \phi \bigr\};
\]
a set of this form is called a \emph{$\Delta$-domain} at $e^{-1}$.

\begin{lemma}[Uniform estimates]\label{lem:uniform}
For every integer $m \ge 0$ there exist a $\Delta$-domain $\mathcal{D}$ at $e^{-1}$ and a constant $C$ such that the following holds.
\begin{enumerate}
\item For all $z \in \mathcal{D}$ and all integers $k, k' \ge 0$ with $|k - k'| \le m$:
\[
\varepsilon_k(z) \ne 0, \qquad |t_k(z)| \le C, \qquad |\varepsilon_k(z)| \le C\, |\varepsilon_{k'}(z)|, \qquad |\Delta_k(z)| \le \frac{C\, |\varepsilon_{k'}(z)|}{k+1} .
\]
\item For all integers $a, b, p, q, r \in \{0, \dots, m\}$ the series
\[
\sum_{k \ge 0} \frac{\Delta_{k+a}(z)\, \Delta_{k+b}(z)}{\varepsilon_{k+p}(z)\, \varepsilon_{k+q}(z)\, \varepsilon_{k+r}(z)}
\]
converges locally uniformly in $\mathcal{D}$, and its sum equals $\frac{1}{4} L(z) + \kappa + o(1)$ as $z \to e^{-1}$ in $\mathcal{D}$, for a constant $\kappa = \kappa(a,b,p,q,r)$.
\end{enumerate}
\end{lemma}

The proof is given in Appendix~\ref{app:uniform}. It is based on the recurrence satisfied by $t - t_h$, in the spirit of the analysis of the height of trees in \cite{flajolet1982}; we give a complete argument because the estimates are needed uniformly in a $\Delta$-domain, for all $h \ge 0$, with index shifts, and (in Appendix~\ref{app:crown}) for a second solution of the same recurrence. We first describe the idea.

\begin{remark}\label{rem:scaling}
In terms of $y = \sqrt{2(1-ez)}$ and $v = hy/2$, the estimates proved below say that, for $z$ close to $e^{-1}$ and $h$ large, $\varepsilon_h(z) \approx y \coth(v)$ and $\Delta_h(z) \approx y^2 / (2 \sinh^2(v))$, so that the summand in Lemma~\ref{lem:uniform}(2) is approximately $\frac{y}{4}\, g_3(v)$ with
\begin{equation}
\label{g3v}
g_3(v) := \frac{1}{\sinh(v) \cosh^3(v)} .
\end{equation}
For $1 \ll h \ll 1/|y|$ this is close to $\frac{1}{2h}$, which is the source of the logarithmic singularity.
\end{remark}

\paragraph{Outline of the proof.} Put $d_h := t - t_h$ and $\tau := 1 - t$, so that $\varepsilon_h = \tau + d_h$. Since $t_{h+1} = z e^{t_h} = t e^{-d_h}$, we have $d_{h+1} = t (1 - e^{-d_h})$. Expanding $1/(1 - e^{-x}) = 1/x + 1/2 + \mathcal{O}(x)$ and iterating, we get
\[
\frac{t^h}{d_h} = \frac{1 - t^h}{2\tau} + b_h ,
\]
where $b_h$ consists of the initial term $1/d_0$ and of a sum of corrections of size $\mathcal{O}(|t^j d_j|)$, $j < h$. The proof has three steps. Away from $e^{-1}$ we have $|t| \le \theta < 1$, and $d_h$ decays geometrically. Near $e^{-1}$, in a $\Delta$-domain, we have $t = e^{-s}$ with $s$ in a sector $|\arg s| \le \pi/2 - \phi_0$, and an induction on $h$ shows that $b_h$ grows at most logarithmically in $h$. Hence $d_h$ is close to $2\tau t^h / (1 - t^h)$ and $|\varepsilon_h|$ is of order $\max(1/h, |s|)$, uniformly in $h$ and $z$; this is the scaling form of Remark~\ref{rem:scaling}, since $s \sim \tau \sim y$, and it gives part~(1). For part~(2), every term of the series is compared with the term obtained from this approximation, which is $4 \tau\, t^{2k} / \bigl( (1 - t^k)(1 + t^k)^3 \bigr)$. The differences are bounded by a summable sequence that does not depend on $z$, so their sum tends to a constant, and the sum of the approximating terms is evaluated explicitly: it equals $\frac{1}{2} \ln(1/\tau)$ plus a constant plus $o(1)$, that is, $\frac{1}{4} L(z)$ plus a constant plus $o(1)$.

\subsection{Proof of Theorem~\ref{th_main}}

By \eqref{eq:ties_discrete},
\[
M_{\mathrm{ties}}(z) = \sum_{h \ge 1} \frac{\Delta_h^2}{\varepsilon_{h-1}^2\, \varepsilon_h} = \sum_{k \ge 0} \frac{\Delta_{k+1}^2}{\varepsilon_k^2\, \varepsilon_{k+1}},
\]
which is the series of Lemma~\ref{lem:uniform}(2) with $a = b = r = 1$ and $p = q = 0$. Let $\mathcal{D}$ be the $\Delta$-domain of that lemma for $m = 1$. The series converges locally uniformly in $\mathcal{D}$, which contains the disc $|z| < e^{-1}$; by Weierstrass' theorem its sum is analytic in $\mathcal{D}$ and its Taylor coefficients at $0$ are obtained termwise, so it is the analytic continuation of the power series $M_{\mathrm{ties}}(z)$. Moreover $M_{\mathrm{ties}}(z) = \frac{1}{4} L(z) + \kappa + o(1)$ as $z \to e^{-1}$ in $\mathcal{D}$, and Theorem~\ref{thm:transfer} yields
\[
[z^n] M_{\mathrm{ties}}(z) = \frac{e^n}{4n}\, (1 + o(1)) .
\]
With $[z^n] M(z) = n^n / n! \sim e^n / \sqrt{2\pi n}$, \eqref{eq:tieprob} gives
\[
\mathbb{P}(\text{$g$ has no unique highest tree}) = \frac{\sqrt{2\pi}}{4}\, n^{-1/2}\, (1 + o(1)) = \sqrt{\frac{\pi}{8}}\; n^{-1/2}\, (1 + o(1)),
\]
which proves Theorem~\ref{th_main}. \qed

\begin{remark}
The proof gives the error term $o(n^{-1/2})$. An explicit error term would follow from a quantitative version of the last step of the proof in Appendix~\ref{app:uniform}, that is, from a rate of convergence in \eqref{eq:kappa1}; we do not pursue this here.
\end{remark}

\section{Generalization to $c$-Branches and Crown Size}
\label{sec:branches}

Throughout this section $c \ge 0$ is a fixed integer.

\subsection{Branches, gaps and exact generating functions}

For a mapping $g \in \Sigma_n$, the \emph{depth} of a vertex is its distance to the cycle core; the cyclic vertices are exactly the vertices of depth $0$.

\begin{definition}\label{def:branch}
A \emph{$c$-branch} of $g$ is the subtree consisting of a vertex $x$ of depth $c$ (its \emph{root}) together with all vertices whose path to the cycle core passes through $x$. Its \emph{height} is the largest distance from $x$ to one of its vertices. In particular, the $0$-branches are the trees of $g$.
\end{definition}

\begin{definition}\label{def:gap}
The \emph{$c$-gap} of $g$ is
\[
\operatorname{gap}_c(g) =
\begin{cases}
0 & \text{no $c$-branch, or the maximal height is attained at least twice,}\\
H_1 - H_2 & \text{a unique highest $c$-branch, of height $H_1$; second largest height $H_2$,}\\
\infty & \text{exactly one $c$-branch.}
\end{cases}
\]
If $\operatorname{gap}_c(g) \ge 1$, we write $T_c$ for the unique highest $c$-branch and $h_{2,c}(g)$ for the largest height among the other $c$-branches, with $h_{2,c}(g) = -1$ if there are none.
\end{definition}

Thus, for $j \ge 1$, $\operatorname{gap}_c(g) \ge j$ means that $g$ has a $c$-branch whose height exceeds the height of every other $c$-branch by at least $j$, and $\operatorname{gap}_0(g) \ge 1$ means that $g$ has a unique highest tree.

\paragraph{Notation.} Recall from Section~\ref{sec:uniform} the conventions $t_{-1} = 0$ and $\Delta_0 = t_0 = z$, and the notation $\varepsilon_h = 1 - t_h$ for $h \ge -1$.
Let $\Phi(x) := z e^{x}$, regarded as a function of $x$ with parameter $z$, and let $\Phi^{(i)}$ be its $i$-fold iterate, $\Phi^{(0)}(x) = x$. Since $t_h = \Phi(t_{h-1})$ for all $h \ge 0$,
\begin{equation}\label{eq:iter}
\Phi^{(i)}(t_h) = t_{h+i} \qquad (h \ge -1,\ i \ge 0).
\end{equation}
By the chain rule and $\Phi' = \Phi$, the derivative of $\Phi^{(c)}$ with respect to $x$ is
\begin{equation}\label{eq:iterprime}
\bigl(\Phi^{(c)}\bigr)'(x) = \prod_{i=1}^{c} \Phi^{(i)}(x).
\end{equation}

\begin{lemma}[Branch decomposition]\label{lem:branch}
Let $\mathcal{X}$ and $\mathcal{Y}$ be disjoint classes of rooted labelled trees with EGFs $X(z)$ and $Y(z)$.
\begin{enumerate}
\item The EGF of mappings all of whose $c$-branches belong to $\mathcal{X}$ is $\bigl(1 - \Phi^{(c)}(X(z))\bigr)^{-1}$.
\item The EGF of mappings having exactly one $c$-branch in $\mathcal{Y}$ and all other $c$-branches in $\mathcal{X}$ is
\[
Y(z)\, \frac{\bigl(\Phi^{(c)}\bigr)'(X(z))}{\bigl(1 - \Phi^{(c)}(X(z))\bigr)^{2}} .
\]
\end{enumerate}
\end{lemma}

\begin{proof}
We first show by induction on $c$ that the EGF of trees all of whose $c$-branches belong to $\mathcal{X}$ is $\Phi^{(c)}(X(z))$. For $c = 0$ a tree is its own unique $0$-branch. For the induction step, a tree is a root with a set of subtrees attached, and the $(c+1)$-branches of the tree are exactly the $c$-branches of these subtrees; hence the EGF is $z \exp\bigl(\Phi^{(c)}(X(z))\bigr) = \Phi^{(c+1)}(X(z))$. Since a mapping is a set of cycles of trees and its $c$-branches are the $c$-branches of its trees, the first claim follows as in the derivation of $M(z)$.

For the second claim, apply the first claim to the class $\mathcal{X} \cup \mathcal{Y}$ with an additional variable $a$ marking the $c$-branches that belong to $\mathcal{Y}$; as the classes are disjoint, this gives $\bigl(1 - \Phi^{(c)}(X(z) + a Y(z))\bigr)^{-1}$. The coefficient of $a^1$ is the stated expression.
\end{proof}

For $k \ge -1$ define
\begin{equation}\label{eq:psi}
\pi_{c,k}(z) := \prod_{i=1}^{c} t_{k+i}(z), \qquad \psi_c(k) := \frac{\pi_{c,k}(z)}{\varepsilon_{k+c}(z)^{2}},
\end{equation}
where the empty product ($c = 0$) equals $1$.

\begin{corollary}\label{cor:branch}
Let $k \ge -1$ and $H > k$.
\begin{enumerate}
\item The EGF of mappings all of whose $c$-branches have height at most $k$ is $1/\varepsilon_{k+c}(z)$. For $k = -1$ these are the mappings without $c$-branches.
\item The EGF of mappings having one $c$-branch of height exactly $H$ and all other $c$-branches of height at most $k$ is $\Delta_H(z)\, \psi_c(k)$.
\end{enumerate}
\end{corollary}

\begin{proof}
Apply Lemma~\ref{lem:branch} with $\mathcal{X}$ the class of trees of height at most $k$ (empty for $k = -1$), whose EGF is $t_k$, and $\mathcal{Y}$ the class of trees of height exactly $H$, whose EGF is $\Delta_H$. By \eqref{eq:iter} and \eqref{eq:iterprime}, $\Phi^{(c)}(t_k) = t_{k+c}$ and $\bigl(\Phi^{(c)}\bigr)'(t_k) = \pi_{c,k}$.
\end{proof}

For $i \ge 0$ let $E^{(c)}_i(z)$ denote the EGF of the mappings $g$ with $\operatorname{gap}_c(g) = i$.

\begin{lemma}[Exact gap generating functions]\label{lem:gap_exact}
As formal power series,
\begin{align}
E^{(c)}_0(z) &= \frac{1}{\varepsilon_{c-1}} + \sum_{k \ge 0} \left[ \frac{\Delta_k\, \pi_{c,k-1}\, \Delta_{k+c}}{\varepsilon_{k+c}\, \varepsilon_{k+c-1}^{2}} + \frac{\rho_{c,k}}{\varepsilon_{k+c}\, \varepsilon_{k+c-1}} \right], \label{eq:E0}\\
E^{(c)}_i(z) &= \sum_{k \ge 0} \Delta_{k+i} \bigl[ \psi_c(k) - \psi_c(k-1) \bigr] \qquad (i \ge 1), \label{eq:Ei}
\end{align}
where $\rho_{c,k} := \Delta_{k+c} - \Delta_k\, \pi_{c,k-1}$, and
\begin{equation}\label{eq:psidiff}
\psi_c(k) - \psi_c(k-1) = \pi_{c,k}\, \Delta_{k+c} \left[ \frac{1}{\varepsilon_{k+c}\, \varepsilon_{k+c-1}^{2}} + \frac{1}{\varepsilon_{k+c}^{2}\, \varepsilon_{k+c-1}} \right] + \frac{t_{k+c} - t_k}{\varepsilon_{k+c-1}^{2}} \prod_{i=1}^{c-1} t_{k+i} .
\end{equation}
Moreover,
\begin{equation}\label{eq:rho}
\rho_{c,k} = \Delta_k^{2} \int_0^1 (1-s)\, \bigl(\Phi^{(c)}\bigr)''\bigl(t_{k-1} + s \Delta_k\bigr)\, ds ,
\end{equation}
so that $\rho_{0,k} = 0$.
\end{lemma}

For $c = 0$ the last term of \eqref{eq:psidiff} vanishes, $\pi_{0,k} = 1$, and \eqref{eq:E0} reduces to $1 + \sum_{k \ge 0} \Delta_k^{2} / (\varepsilon_k\, \varepsilon_{k-1}^{2})$, in agreement with \eqref{eq:ties_discrete}.

\begin{proof}
Let $i \ge 1$. A mapping has $\operatorname{gap}_c = i$ if and only if, for some $k \ge 0$, it has one $c$-branch of height exactly $k + i$, all other $c$-branches have height at most $k$, and not all of them have height at most $k - 1$. By Corollary~\ref{cor:branch}(2) the EGF of these mappings is $\Delta_{k+i} [\psi_c(k) - \psi_c(k-1)]$, which gives \eqref{eq:Ei}. Writing
\[
\psi_c(k) - \psi_c(k-1) = \pi_{c,k} \left[ \frac{1}{\varepsilon_{k+c}^{2}} - \frac{1}{\varepsilon_{k+c-1}^{2}} \right] + \frac{\pi_{c,k} - \pi_{c,k-1}}{\varepsilon_{k+c-1}^{2}}
\]
and using $\varepsilon_{k+c-1} - \varepsilon_{k+c} = \Delta_{k+c}$ and $\pi_{c,k} - \pi_{c,k-1} = (t_{k+c} - t_k) \prod_{i=1}^{c-1} t_{k+i}$, we obtain \eqref{eq:psidiff}.

A mapping has $\operatorname{gap}_c = 0$ if and only if it has no $c$-branch, or for some $k \ge 0$ the maximal height of its $c$-branches equals $k$ and is attained at least twice. By Corollary~\ref{cor:branch}, the EGF of the mappings of the second kind for a given $k$ is
\[
\frac{1}{\varepsilon_{k+c}} - \frac{1}{\varepsilon_{k+c-1}} - \Delta_k\, \psi_c(k-1)
= \frac{\Delta_{k+c}}{\varepsilon_{k+c}\, \varepsilon_{k+c-1}} - \frac{\Delta_k\, \pi_{c,k-1}}{\varepsilon_{k+c-1}^{2}} .
\]
Substituting $\Delta_{k+c} = \Delta_k \pi_{c,k-1} + \rho_{c,k}$ in the first fraction and using $\varepsilon_{k+c-1} - \varepsilon_{k+c} = \Delta_{k+c}$ once more gives the summand of \eqref{eq:E0}. Finally, by \eqref{eq:iter} and \eqref{eq:iterprime},
\[
\rho_{c,k} = \Phi^{(c)}(t_{k-1} + \Delta_k) - \Phi^{(c)}(t_{k-1}) - \Delta_k \bigl(\Phi^{(c)}\bigr)'(t_{k-1}),
\]
and \eqref{eq:rho} is Taylor's formula with integral remainder.
\end{proof}

\subsection{Asymptotics of the gap}

The only analytic input needed in this section is Lemma~\ref{lem:uniform}; recall that $L(z) = \ln \frac{1}{1-ez}$.

\begin{proposition}\label{prop:gap_sing}
Let $i \ge 1$. There is a $\Delta$-domain $\mathcal{D}$ at $e^{-1}$ in which $E^{(c)}_0$ and $E^{(c)}_i$ are analytic, and, as $z \to e^{-1}$ in $\mathcal{D}$,
\[
E^{(c)}_0(z) = \tfrac{1}{4} L(z) + \kappa_0 + o(1), \qquad E^{(c)}_i(z) = \tfrac{1}{2} L(z) + \kappa_i + o(1),
\]
for some constants $\kappa_0, \kappa_i$ depending on $c$ and $i$.
\end{proposition}

\begin{proof}
Let $\mathcal{D}$ and $C$ be as in Lemma~\ref{lem:uniform} with $m = c + i + 1$. Call a series $\sum_{k \ge 1} f_k$ of functions analytic in $\mathcal{D}$ \emph{tame} if every $f_k$ extends continuously to $z = e^{-1}$ and $|f_k(z)| \le C' k^{-2}$ on $\mathcal{D}$ for a constant $C'$. A tame series converges uniformly on $\mathcal{D} \cup \{e^{-1}\}$; hence its sum is analytic in $\mathcal{D}$ and equals $\kappa + o(1)$ as $z \to e^{-1}$ in $\mathcal{D}$, where $\kappa$ is its value at $e^{-1}$.

All terms occurring below are rational expressions in the entire functions $t_h$ and $\rho_{c,k}$ whose denominators are products of functions $\varepsilon_h$; these do not vanish in $\mathcal{D}$ by Lemma~\ref{lem:uniform}(1), nor at $e^{-1}$, because $t_h(e^{-1}) < t(e^{-1}) = 1$. Hence every single term is analytic in $\mathcal{D}$ and continuous at $e^{-1}$; in particular we may set aside the terms with $k = 0$ and the term $1/\varepsilon_{c-1}$. By Lemma~\ref{lem:uniform}(1), for bounded analytic weights $w_k$ and shifts $a, b, p, q \in \{0, \dots, m\}$,
\begin{equation}\label{eq:tame}
\sum_{k \ge 1} w_k\, \frac{\Delta_{k+a}\, \Delta_{k+b}}{\varepsilon_{k+p}\, \varepsilon_{k+q}} \quad \text{is tame,}
\end{equation}
since each factor $\Delta_{k+a} / \varepsilon_{k+p}$ is bounded by $C/(k+1)$.

\emph{The case $i \ge 1$.} Insert \eqref{eq:psidiff} into \eqref{eq:Ei}. The contribution of the last term of \eqref{eq:psidiff} is tame by \eqref{eq:tame}, because $t_{k+c} - t_k = \sum_{l=1}^{c} \Delta_{k+l}$ and $|t_h| \le C$. In the remaining two terms write $\pi_{c,k} = 1 + (\pi_{c,k} - 1)$. Since $|t_h| \le C$, we have $|\pi_{c,k} - 1| \le C'' \sum_{l=1}^{c} |\varepsilon_{k+l}|$, and $|\varepsilon_{k+l}| \le C |\varepsilon_{k+c-1}|$, so the contribution of $\pi_{c,k} - 1$ is again tame by \eqref{eq:tame}. What remains is
\[
\sum_{k \ge 1} \frac{\Delta_{k+i}\, \Delta_{k+c}}{\varepsilon_{k+c}\, \varepsilon_{k+c-1}^{2}} + \sum_{k \ge 1} \frac{\Delta_{k+i}\, \Delta_{k+c}}{\varepsilon_{k+c}^{2}\, \varepsilon_{k+c-1}} ,
\]
and each of the two sums equals $\frac{1}{4} L(z) + \kappa + o(1)$ by Lemma~\ref{lem:uniform}(2) (for $c = 0$ after shifting the index of summation by one; the finitely many terms by which the sums then differ are continuous at $e^{-1}$).

\emph{The case $i = 0$.} In \eqref{eq:E0}, the first summand is treated in the same way: replacing $\pi_{c,k-1}$ by $1$ changes the sum by a tame series, and the resulting sum equals $\frac{1}{4} L(z) + \kappa + o(1)$ by Lemma~\ref{lem:uniform}(2). For the second summand, note that $t_{k-1} + s \Delta_k = (1-s) t_{k-1} + s t_k$ has modulus at most $C$ for $0 \le s \le 1$, and that $(\Phi^{(c)})''(x)$ is an entire function of $x$ and $z$, hence bounded for $|x| \le C$ and $z \in \mathcal{D}$. By \eqref{eq:rho}, $|\rho_{c,k}| \le C''' |\Delta_k|^{2}$, and the series $\sum_k \rho_{c,k} / (\varepsilon_{k+c}\, \varepsilon_{k+c-1})$ is tame.

Finally, the series obtained in this way converge locally uniformly in $\mathcal{D}$, which contains the disc $|z| < e^{-1}$; by Weierstrass' theorem the Taylor coefficients of their sums at $0$ are obtained termwise, so the sums are the analytic continuations of the formal power series $E^{(c)}_0$ and $E^{(c)}_i$.
\end{proof}

\begin{theorem}\label{th_gap}
Let $c \ge 0$ and $j \ge 1$ be fixed integers and let $g \in \Sigma_n$ be a random mapping. Then, as $n \to \infty$,
\[
\mathbb{P}\bigl(\operatorname{gap}_c(g) = 0\bigr) \sim \sqrt{\frac{\pi}{8}}\; n^{-1/2}, \qquad
\mathbb{P}\bigl(\operatorname{gap}_c(g) = i\bigr) \sim \sqrt{\frac{\pi}{2}}\; n^{-1/2} \quad (i \ge 1),
\]
and consequently
\[
\mathbb{P}\bigl(\operatorname{gap}_c(g) < j\bigr) \sim (2j - 1) \sqrt{\frac{\pi}{8}}\; n^{-1/2} .
\]
In particular, with probability $1 - \Theta(n^{-1/2})$ the highest $c$-branch of $g$ is unique and exceeds every other $c$-branch in height by at least $j$.
\end{theorem}

\begin{proof}
We have $[z^n] L(z) = e^n / n$ for $n \ge 1$, constants do not contribute to $[z^n]$ for $n \ge 1$, and a function that is analytic in a $\Delta$-domain at $e^{-1}$ and is $o(1)$ as $z \to e^{-1}$ has coefficients $o(e^n / n)$ by the transfer theorem \cite{flajolet1990}. Hence Proposition~\ref{prop:gap_sing} yields
\[
[z^n] E^{(c)}_0(z) \sim \frac{e^n}{4n}, \qquad [z^n] E^{(c)}_i(z) \sim \frac{e^n}{2n} \quad (i \ge 1).
\]
Dividing by $[z^n] M(z) = n^n / n! \sim e^n / \sqrt{2 \pi n}$ gives the first two claims, and the third follows by summing over $0 \le i < j$.
\end{proof}

The constants in Theorem~\ref{th_gap} do not depend on $c$: asymptotically, each positive value of the gap is twice as likely as a tie. For $c = 0$ and $j = 1$ the theorem recovers the leading term of the result of Section~\ref{sec:unique_tree}.

\subsection{Crown size}
\label{sec:crown}

\begin{definition}
Let $\operatorname{gap}_c(g) \ge 1$. The \emph{$c$-crown} of $g$, denoted $H$, is the forest consisting of all vertices of $T_c$ at distance at least $h_{2,c}(g) + 1$ from the root of $T_c$, and $r$ is the number of its roots, that is, of the vertices of $T_c$ at distance exactly $h_{2,c}(g) + 1$ from the root of $T_c$. If $\operatorname{gap}_c(g) = 0$, we put $H = \emptyset$ and $r = 0$.
\end{definition}

\begin{theorem}
\label{th_high_tree}
Let $c \ge 0$ be a fixed integer and $\alpha > 1$ a fixed real number. Let $g \in \Sigma_n$ be a random mapping and let $H$ be the $c$-crown of $g$, having $r$ roots. Then $|H| > \alpha r > 0$ with probability $1 - \Theta(n^{-1/2})$.
\end{theorem}

Throughout this subsection $c \ge 0$ and $\alpha > 1$ are fixed. The complement of the event $|H| > \alpha r > 0$ is the disjoint union of the event $\operatorname{gap}_c(g) = 0$ (equivalently, $r = 0$) and the event
\[
\mathcal{S} := \bigl\{ \operatorname{gap}_c(g) \ge 1 \text{ and } |H| \le \alpha r \bigr\}.
\]
By Theorem~\ref{th_gap}, the first event has probability $\Theta(n^{-1/2})$, so the task is to show that $\mathbb{P}(\mathcal{S}) = \mathcal{O}(n^{-1/2})$. Let $M_{\mathrm{sparse}}(z)$ denote the EGF of the mappings in $\mathcal{S}$. For power series $F$ and $G$ we write $F \preceq G$ if $[z^n] F \le [z^n] G$ for all $n$. We bound $M_{\mathrm{sparse}}$ from above by a series that can be analysed with the tools of Section~\ref{sec:unique_tree}. Put
\[
\Psi_k(z) := \psi_c(k) - \psi_c(k-1) \quad (k \ge 0), \qquad p_\alpha(z) := \sum_{1 \le \nu \le \alpha} \frac{\nu^{\nu-1}}{\nu!} z^\nu ,
\]
see \eqref{eq:psi} and \eqref{eq:psidiff}.
By Theorem~\ref{thm:tree}, $\lambda^\alpha\, t\bigl(1/(e\lambda)\bigr) = 1 - \sqrt{2(\lambda - 1)}\,(1 + o(1))$ as $\lambda \downarrow 1$. We fix $\lambda > 1$ such that
\begin{equation}\label{eq:S0}
S_0 := \lambda^\alpha\, t\bigl(1/(e\lambda)\bigr) \in (0, 1), \qquad \text{and put} \qquad S(z) := \lambda^\alpha\, t(z/\lambda),
\end{equation}
a function analytic for $|z| < \lambda e^{-1}$ with non-negative Taylor coefficients. Let
\[
\tilde t_k(z) := \Phi^{(k)}\bigl( z e^{S(z)} \bigr), \qquad \hat Y_k(z) := \tilde t_k(z) - t_k(z) \qquad (k \ge 0).
\]
Finally let
\[
\hat M(z) := \sum_{k \ge 0} \hat Y_k(z)\, \Psi_k(z) .
\]

\begin{lemma}[Majorant]\label{lem:majorant}
We have $\hat Y_k \succeq 0$ and $\Psi_k \succeq 0$ for all $k \ge 0$, and
\begin{equation}\label{eq:maj}
M_{\mathrm{sparse}}(z) \preceq p_\alpha(z)\, \psi_c(-1) + \hat M(z) .
\end{equation}
\end{lemma}

\begin{proof}
For $k \ge 0$ let $\mathcal{Y}_k$ be the class of trees of height at least $k+1$ in which the number $N$ of vertices at distance at least $k+1$ from the root and the number $r$ of vertices at distance exactly $k+1$ from the root satisfy $N \le \alpha r$, and let $Y_k(z)$ be its EGF.
A mapping $g$ with $h_{2,c}(g) = k \ge 0$ belongs to $\mathcal{S}$ if and only if it has one $c$-branch in $\mathcal{Y}_k$, all its other $c$-branches have height at most $k$, and not all of them have height at most $k-1$. A mapping $g$ with $h_{2,c}(g) = -1$ belongs to $\mathcal{S}$ if and only if it has exactly one $c$-branch and this $c$-branch has at most $\alpha$ vertices. Hence, by Lemma~\ref{lem:branch}(2) and \eqref{eq:iter}, \eqref{eq:iterprime},
\begin{equation}\label{eq:sparse_exact}
M_{\mathrm{sparse}}(z) = p_\alpha(z)\, \psi_c(-1) + \sum_{k \ge 0} Y_k(z)\, \Psi_k(z) .
\end{equation}
Note that $\Psi_k \succeq 0$, because $\psi_c(k)$ is a power series with non-negative coefficients in $z$ and $t_k$, and $t_{k-1} \preceq t_k$.

We claim that
\begin{equation}\label{eq:Yk}
0 \preceq Y_k \preceq \hat Y_k \qquad (k \ge 0).
\end{equation}
Indeed, the EGF of the trees of height at most $k$, with an additional variable $u$ marking the vertices at distance $k$ from the root, is $\Phi^{(k)}(zu)$ (by induction on $k$, as in the proof of Lemma~\ref{lem:branch}). Every tree is obtained in a unique way from a tree of height at most $k$ by attaching a set of trees to each of its vertices at distance $k$ from the root; the resulting tree has height at most $k$ if and only if all these sets are empty. Give every attached tree with $\nu$ vertices the weight $\lambda^{\alpha - \nu}$, and every tree the product of the weights of the trees attached in this way. The weighted EGF of a single attached tree is $\sum_\nu \lambda^{\alpha-\nu} \frac{\nu^{\nu-1}}{\nu!} z^\nu = S(z)$, so the weighted EGF of all trees is $\Phi^{(k)}(z e^{S(z)}) = \tilde t_k$, and that of the trees of height at least $k+1$ is $\tilde t_k - t_k = \hat Y_k$. A tree of height at least $k+1$ has weight $\lambda^{\alpha r - N}$, which is positive, and at least $1$ if the tree belongs to $\mathcal{Y}_k$. This proves \eqref{eq:Yk}.

Now \eqref{eq:maj} follows from \eqref{eq:sparse_exact}, \eqref{eq:Yk} and $\Psi_k \succeq 0$.
\end{proof}

The weights $\lambda^{\alpha r - N}$ used in this proof are those of a Chernoff bound for the event $N \le \alpha r$. The reason why $\hat M$ can be analysed is that the functions $\tilde t_k$ satisfy the same recurrence $\tilde t_{k+1} = z e^{\tilde t_k}$ as the functions $t_k$, with the initial value $z e^{S(z)}$ instead of $z$; the condition $S_0 < 1$ in \eqref{eq:S0} ensures that this initial value is smaller than $t(e^{-1}) = 1$ at $z = e^{-1}$. The estimates of Appendix~\ref{app:uniform} therefore apply to $t - \tilde t_k$ as well, and $\hat Y_k = \tilde t_k - t_k$ turns out to be asymptotically proportional to $\Delta_{k+1}$. This leads to the following statement, which is proved in Appendix~\ref{app:crown}.

\begin{proposition}\label{prop:Mhat}
There exist a $\Delta$-domain $\mathcal{D}$ at $e^{-1}$ and a constant $\beta^*$ such that the series $\hat M(z)$ converges locally uniformly in $\mathcal{D}$ and $\hat M(z) - \beta^* L(z)$ is bounded in $\mathcal{D}$.
\end{proposition}

\begin{proof}[Proof of Theorem~\ref{th_high_tree}]
The lower bound for the probability of the complement of the event $|H| > \alpha r > 0$ follows from Theorem~\ref{th_gap}, as noted above. For the upper bound we use \eqref{eq:maj}. The function $p_\alpha\, \psi_c(-1) = p_\alpha\, \pi_{c,-1} / \varepsilon_{c-1}^2$ is analytic in a disc of radius larger than $e^{-1}$, since $\varepsilon_{c-1}$ is entire and does not vanish for $|z| \le e^{-1}$ (because $|t_{c-1}(z)| \le t_{c-1}(e^{-1}) < 1$ there); hence its coefficients are $\mathcal{O}(\varrho^{n})$ for some $\varrho < e$. By Proposition~\ref{prop:Mhat} and Weierstrass' theorem, the sum of the series $\hat M$ is the analytic continuation of the power series $\hat M(z)$ to $\mathcal{D}$, and by Theorem~\ref{thm:transfer}
\[
[z^n] \hat M(z) = \mathcal{O}\Bigl( \frac{e^n}{n} \Bigr) .
\]
Dividing by $[z^n] M(z) \sim e^n / \sqrt{2\pi n}$ we conclude that $\mathbb{P}(\mathcal{S}) = \mathcal{O}(n^{-1/2})$. Hence the complement of the event $|H| > \alpha r > 0$ has probability $\Theta(n^{-1/2})$.
\end{proof}

\section{Numerical Verification}
\label{sec:empitrical}

We checked the results in two independent ways: by evaluating the exact generating functions of Section~\ref{sec:branches}, and by simulating random mappings.
The programs and the raw results are available at \url{https://github.com/berlm/unique-highest-tree}.

\subsection{Exact values}

By Corollary~\ref{cor:branch}, the EGF of the mappings with $\operatorname{gap}_c \ge j$ is $\sum_{H \ge 0} \Delta_H\, \psi_c(\max(H-j,-1))$. We computed the coefficients of these series for $n \le 25{,}000$ by power series arithmetic in double precision, in the scaled variable $ez$. The sum over $H$ was stopped when the probability that a tree with $25{,}000$ vertices has height larger than $H$ had dropped below $10^{-13}$, so the values are exact up to errors of this order. For $n \le 5$ they agree with a direct enumeration of all $n^n$ mappings.

Let $p^{(c)}_n := \mathbb{P}(\operatorname{gap}_c(g) = 0)$ for a random $g \in \Sigma_n$; thus $1 - p^{(0)}_n$ is the probability of a unique highest tree. Table~\ref{tab:exact_ties} shows $\sqrt{n}\, p^{(c)}_n$ for $c = 0, 1, 2$, which by Theorems~\ref{th_main} and~\ref{th_gap} tends to $\sqrt{\pi/8} \approx 0.62666$ for every $c$. The last three columns show the deviation from the leading term on the scale $1/n$. They suggest that
\[
p^{(c)}_n = \sqrt{\frac{\pi}{8}}\; n^{-1/2} + \frac{c-1}{2n} + o\Bigl(\frac{1}{n}\Bigr) \qquad (c = 0, 1, 2),
\]
which we have not proved. In particular, the convergence in Theorem~\ref{th_main} appears to be of order $1/n$.

\begin{table}[htbp]
\centering
\caption{Exact tie probabilities $p^{(c)}_n = \mathbb{P}(\operatorname{gap}_c = 0)$; here $q_n := \sqrt{\pi/(8n)}$ is the leading term.}
\label{tab:exact_ties}
\begin{tabular}{rcccccc}
\hline
 & \multicolumn{3}{c}{$\sqrt{n}\; p^{(c)}_n$} & \multicolumn{3}{c}{$n\, \bigl(p^{(c)}_n - q_n\bigr)$} \\
$n$ & $c=0$ & $c=1$ & $c=2$ & $c=0$ & $c=1$ & $c=2$ \\ \hline
100 & 0.58004 & 0.62718 & 0.68075 & $-0.4661$ & $+0.0052$ & $+0.5410$ \\
250 & 0.59643 & 0.62687 & 0.65985 & $-0.4780$ & $+0.0033$ & $+0.5248$ \\
1{,}000 & 0.61120 & 0.62671 & 0.64285 & $-0.4887$ & $+0.0017$ & $+0.5120$ \\
2{,}500 & 0.61680 & 0.62668 & 0.63681 & $-0.4928$ & $+0.0010$ & $+0.5075$ \\
10{,}000 & 0.62169 & 0.62666 & 0.63169 & $-0.4964$ & $+0.0005$ & $+0.5037$ \\
25{,}000 & 0.62351 & 0.62666 & 0.62983 & $-0.4978$ & $+0.0002$ & $+0.5022$ \\
\hline
limit & 0.62666 & 0.62666 & 0.62666 & & & \\ \hline
\end{tabular}
\end{table}

Table~\ref{tab:exact_gaps} shows $\sqrt{n}\; \mathbb{P}(\operatorname{gap}_c < j)$ for $j = 2, 3$, which by Theorem~\ref{th_gap} tends to $(2j-1)\sqrt{\pi/8}$. As for $j = 1$, the deviation of the scaled quantity from its limit is of order $n^{-1/2}$, with a larger constant.

\begin{table}[htbp]
\centering
\caption{Exact values of $\sqrt{n}\; \mathbb{P}(\operatorname{gap}_c < j)$.}
\label{tab:exact_gaps}
\begin{tabular}{rcccccc}
\hline
 & \multicolumn{3}{c}{$j = 2$} & \multicolumn{3}{c}{$j = 3$} \\
$n$ & $c=0$ & $c=1$ & $c=2$ & $c=0$ & $c=1$ & $c=2$ \\ \hline
100 & 1.6349 & 1.7559 & 1.8908 & 2.5554 & 2.7267 & 2.9145 \\
250 & 1.7169 & 1.7994 & 1.8882 & 2.7399 & 2.8634 & 2.9953 \\
1{,}000 & 1.7946 & 1.8388 & 1.8847 & 2.9227 & 2.9923 & 3.0645 \\
2{,}500 & 1.8250 & 1.8537 & 1.8831 & 2.9965 & 3.0426 & 3.0899 \\
10{,}000 & 1.8521 & 1.8667 & 1.8816 & 3.0633 & 3.0873 & 3.1116 \\
25{,}000 & 1.8622 & 1.8716 & 1.8810 & 3.0886 & 3.1040 & 3.1196 \\
\hline
limit & 1.8800 & 1.8800 & 1.8800 & 3.1333 & 3.1333 & 3.1333 \\ \hline
\end{tabular}
\end{table}

\subsection{Simulation}

For the crown we have no exact formula, so Theorem~\ref{th_high_tree} was checked by simulation. For each $n$ we generated independent uniformly random mappings and determined, for $c = 0, 1, 2$, the gap and the crown of the highest $c$-branch. Each mapping is processed in linear time: removing leaves repeatedly (Kahn's algorithm) yields the cyclic vertices, the height of the subtree below every vertex, and an order of the vertices from which depths and branch membership are obtained. For $n \le 5$ the program was checked against an independent enumeration of all mappings.

As a check of the simulation, the observed numbers of mappings with $\operatorname{gap}_c < j$ were compared with the exact probabilities of the previous subsection for all six values of $n$, $c \in \{0,1,2\}$ and $j \in \{1,2,3\}$. The largest deviation among these $54$ (dependent) comparisons is $1.76$ standard deviations.

Table~\ref{tab:crown} shows, for $c = 1$ and $\alpha \in \{1.5, 2, 3\}$, the observed frequency of the complement of the event $|H| > \alpha r > 0$, multiplied by $\sqrt{n}$, with $95\%$ Wilson score intervals. By Theorem~\ref{th_high_tree} this quantity is bounded above and below by positive constants. The observed values increase with $n$ and level off, more slowly for larger $\alpha$; this is consistent with the theorem, but the data do not determine the limits. The results for $c = 0$ and $c = 2$ are similar and are included with the published data.

\begin{table}[htbp]
\centering
\caption{Simulated values of $\sqrt{n}\; \mathbb{P}\bigl(\text{not } |H| > \alpha r > 0\bigr)$ for $1$-branches, with half-widths of the $95\%$ confidence intervals.}
\label{tab:crown}
\begin{tabular}{rrccc}
\hline
$n$ & trials & $\alpha = 1.5$ & $\alpha = 2$ & $\alpha = 3$ \\ \hline
100 & 4{,}000{,}000 & 2.208 $\pm$ 0.004 & 3.046 $\pm$ 0.005 & 4.077 $\pm$ 0.005 \\
250 & 4{,}000{,}000 & 2.336 $\pm$ 0.005 & 3.314 $\pm$ 0.006 & 4.650 $\pm$ 0.007 \\
1{,}000 & 2{,}000{,}000 & 2.452 $\pm$ 0.012 & 3.584 $\pm$ 0.014 & 5.292 $\pm$ 0.016 \\
2{,}500 & 1{,}000{,}000 & 2.506 $\pm$ 0.021 & 3.705 $\pm$ 0.026 & 5.593 $\pm$ 0.031 \\
10{,}000 & 400{,}000 & 2.537 $\pm$ 0.049 & 3.775 $\pm$ 0.059 & 5.830 $\pm$ 0.073 \\
25{,}000 & 200{,}000 & 2.477 $\pm$ 0.086 & 3.780 $\pm$ 0.106 & 5.875 $\pm$ 0.131 \\
\hline
\end{tabular}
\end{table}

\section{Conclusion}
We have shown that a random mapping of $n$ elements has a unique highest tree with probability $1 - \sqrt{\frac{\pi}{8}}\, n^{-1/2} + o(n^{-1/2})$, and, more generally, that for all fixed $c \ge 0$ and $j \ge 1$ the highest $c$-branch exceeds every other $c$-branch in height by at least $j$ except with probability asymptotic to $(2j-1)\sqrt{\frac{\pi}{8}}\, n^{-1/2}$. We have also shown that, with probability $1 - \Theta(n^{-1/2})$, the crown of the highest $c$-branch has more than $\alpha$ times as many vertices as roots, for every fixed $\alpha > 1$. The proofs rest on exact generating functions and on estimates for the generating functions of trees of bounded height that are uniform in a $\Delta$-domain (Lemma~\ref{lem:uniform}).

The exact values and simulations of Section~\ref{sec:empitrical} are consistent with these results. In the context of synchronization, Theorem~\ref{th_high_tree} supplies the property of random mappings that is used in the author's proof of Cameron's conjecture \cite{berlinkov2016}.

Several questions remain open. Our proof of Theorem~\ref{th_main} gives the error term $o(n^{-1/2})$. The exact values in Table~\ref{tab:exact_ties} suggest that it is in fact of order $n^{-1}$, with second term $-\frac{1}{2n}$ in the probability of a tie; an error term of this order would follow from a quantitative version of the last step in the proof of Lemma~\ref{lem:uniform}.
For the crown, Theorem~\ref{th_high_tree} determines only the order of the failure probability, and the constant, as a function of $c$ and $\alpha$, is not known. Finally, all parameters $c$, $j$ and $\alpha$ are fixed here; it would be interesting to find the ranges in which the results remain valid when they grow with $n$.

\section*{Acknowledgements}
The author wants to express his deepest gratitude to Mikhail V. Volkov for his continuous support and encouragement prior to and during the work on the initial version of this paper. The author is also grateful for the inspiration and kind words of Avraham Trahtman, who tragically passed away in 2024. During the preparation of this work, the author extensively used generative AI models (Gemini 1.5 Pro and Claude Opus 5.5) to elaborate on the original idea and an internal Google publishing review tool to find any issues. The author reviewed and edited the content to his best abilities and takes full responsibility for the content of the publication.

\appendix

\section{Proof of Lemma~\ref{lem:uniform}}
\label{app:uniform}

We use the notation of Section~\ref{sec:uniform}.

\paragraph{Preliminaries.} We use the following standard facts about the tree function \cite{flajolet2009}; see also Theorem~\ref{thm:tree}.
\begin{itemize}
\item[(T1)] For every $\phi \in (0, \pi/2)$ there is $\eta_0 > 0$ such that $t$ is analytic in $\mathcal{D}(\eta_0, \phi)$, satisfies $t = z e^{t}$ there, and $t(z) = 1 - y + \mathcal{O}(y^2)$ as $z \to e^{-1}$ in $\mathcal{D}(\eta_0, \phi)$, where $y = \sqrt{2(1-ez)}$ with the principal branch of the square root.
\item[(T2)] $t$ extends continuously to the closed disc $|z| \le e^{-1}$ and analytically to a neighbourhood of every point of this disc other than $e^{-1}$, and $|t(z)| < 1$ for $|z| \le e^{-1}$, $z \ne e^{-1}$. (The last assertion holds because the Taylor coefficients of $t$ are positive and $t(e^{-1}) = 1$.)
\end{itemize}
Put
\[
d_h := t - t_h \quad (h \ge -1), \qquad \tau := 1 - t,
\]
so that $\varepsilon_h = \tau + d_h$ and $\Delta_h = d_{h-1} - d_h$ for $h \ge 0$. Since $t_{h+1} = z e^{t_h} = z e^{t} e^{-d_h} = t e^{-d_h}$,
\begin{equation}\label{eq:rec}
d_{h+1} = t \bigl( 1 - e^{-d_h} \bigr) \qquad (h \ge -1).
\end{equation}
We shall use the functions
\[
R(x) := \frac{1}{1 - e^{-x}} - \frac{1}{x} - \frac{1}{2}, \qquad G(x) := x - 1 + e^{-x} .
\]
$R$ is analytic for $|x| < 2\pi$ with $R(0) = 0$, and $R(x) > 0$ for real $x > 0$. Hence there is a constant $c_R$ with $|R(x)| \le c_R |x|$ for $|x| \le 1$. Moreover, $|G(x)| \le |x|^2$ and $|G(x) - x^2/2| \le |x|^3$ for $|x| \le 1$.

The power series $d_h$ has non-negative coefficients, because $t_h$ counts a subclass of the trees. Hence, with $x_h := d_h(e^{-1}) = 1 - t_h(e^{-1})$,
\begin{equation}\label{eq:inside}
|d_h(z)| \le x_h \quad \text{and} \quad |\varepsilon_h(z)| \ge 1 - t_h(e^{-1}) = x_h \qquad (|z| \le e^{-1},\ h \ge 0).
\end{equation}
By \eqref{eq:rec} at $z = e^{-1}$, $x_0 = 1 - e^{-1}$ and $1/x_{h+1} = 1/x_h + 1/2 + R(x_h) > 1/x_h + 1/2$, so
\begin{equation}\label{eq:xh}
0 < x_h < \frac{2}{h+3} \qquad (h \ge 0).
\end{equation}

\paragraph{Admissible sequences.} For later use in Appendix~\ref{app:crown}, part of the estimates is proved for arbitrary solutions of the recurrence \eqref{eq:rec}. A sequence $(d^*_h)_{h \ge 0}$ of functions is called \emph{admissible} if
\begin{itemize}
\item[(a)] $d^*_0$ is analytic in the domains $\mathcal{D}(\eta_0, \phi)$ of (T1) and in a neighbourhood of every point of the disc $|z| \le e^{-1}$ other than $e^{-1}$;
\item[(b)] $d^*_{h+1} = t \bigl( 1 - e^{-d^*_h} \bigr)$ for all $h \ge 0$;
\item[(c)] $\sup \{ |d^*_h(z)| : |z| \le e^{-1},\ z \ne e^{-1} \} \to 0$ as $h \to \infty$;
\item[(d)] $d^*_0(z)$ tends to a limit $x^*_0 \in (0,1)$ as $z \to e^{-1}$ in $\mathcal{D}(\eta_0, \phi)$.
\end{itemize}
By (b), (d) and $t \to 1$, for every fixed $h$ the function $d^*_h(z)$ tends to a limit $x^*_h$ as $z \to e^{-1}$ in $\mathcal{D}(\eta_0, \phi)$, where $x^*_{h+1} = 1 - e^{-x^*_h}$. As in \eqref{eq:xh}, $1/x^*_h > 1/x^*_0 + h/2 > (h+2)/2$, so
\begin{equation}\label{eq:xstar}
0 < x^*_h < \min\Bigl( 1, \frac{2}{h+2} \Bigr) \qquad (h \ge 0).
\end{equation}
The sequence $(d_h)_{h \ge 0}$ is admissible with $x^*_h = x_h$, by (T1), (T2), \eqref{eq:rec}, \eqref{eq:inside} and \eqref{eq:xh}.

\paragraph{Step 1: away from $e^{-1}$.}
\emph{Claim A. For every $\delta > 0$ and every admissible sequence $(d^*_h)$ there exist $\eta > 0$, $\theta \in (0,1)$ and $C$ such that for all $z$ in $\mathcal{K} := \{ z : |z| \le (1+\eta) e^{-1},\ |z - e^{-1}| \ge \delta \}$ and all $h \ge 0$:}
\[
|d^*_h(z)| \le C \theta^h, \qquad |d_h(z)| \le C \theta^h, \qquad C^{-1} \le |\varepsilon_h(z)| \le C .
\]
Let $\mathcal{K}_0 := \{ z : |z| \le e^{-1},\ |z - e^{-1}| \ge \delta/2 \}$. By (T2), $\theta_0 := \max_{\mathcal{K}_0} |t| < 1$. Fix $\sigma \in (0, (1-\theta_0)/8]$ such that $\theta := (\theta_0 + \sigma) e^{2\sigma} < 1$, and, using \eqref{eq:xh} and (c), fix $h_0$ such that $x_{h_0} \le \sigma$ and $|d^*_{h_0}| \le \sigma$ on $\mathcal{K}_0$. On $\mathcal{K}_0$ we then have $|t| \le \theta_0$, $|d_{h_0}| \le \sigma$, $|d^*_{h_0}| \le \sigma$, and $|\varepsilon_h| \ge x_{h_0}$ for $h \le h_0$, by \eqref{eq:inside}. The functions $t$, $d_{h_0}$, $d^*_{h_0}$ and $\varepsilon_0, \dots, \varepsilon_{h_0}$ are analytic in a neighbourhood of $\mathcal{K}_0$, and every point of $\mathcal{K}$ has distance at most $\eta e^{-1}$ from $\mathcal{K}_0$ as soon as $\eta \le e\delta/2$. By uniform continuity we can therefore choose $\eta$ so small that on $\mathcal{K}$
\[
|t| \le \theta_0 + \sigma, \qquad |d_{h_0}| \le 2\sigma, \qquad |d^*_{h_0}| \le 2\sigma, \qquad |\varepsilon_h| \ge x_{h_0}/2 \quad (h \le h_0).
\]
Since $|1 - e^{-x}| \le |x| e^{|x|}$, the recurrence and induction give $|d_{h+1}| \le \theta |d_h| \le 2\sigma$ for all $h \ge h_0$ on $\mathcal{K}$, hence $|d_h| \le 2\sigma \theta^{h - h_0}$; the same holds for $d^*_h$. Finally $|\tau| \ge 1 - |t| \ge 1 - \theta_0 - \sigma$, so $|\varepsilon_h| \ge |\tau| - |d_h| \ge (1-\theta_0)/2$ for $h \ge h_0$. The upper bounds are clear, since $t$ and the finitely many functions $d_h$, $d^*_h$ with $h \le h_0$ are bounded on $\mathcal{K}$. This proves Claim~A.

\paragraph{Step 2: near $e^{-1}$.}
Fix $\phi \in (0, \pi/2)$, put $\phi_0 := \phi/4$, $\gamma := \sin \phi_0$, and let
\[
\mathcal{W} := \bigl\{ w \in \mathbb{C} \setminus \{0\} : |\arg w| \le \pi/2 - \phi_0 \bigr\},
\]
so that $\operatorname{Re} w \ge \gamma |w|$ for $w \in \mathcal{W}$. There is a constant $\gamma' \in (0, 1]$, depending only on $\phi$, such that for all $w \in \mathcal{W}$
\begin{equation}\label{eq:sector}
\gamma' \min(1, |w|) \le |1 - e^{-w}| \le \min(2, |w|), \qquad |1 + e^{-w}| \ge \gamma', \qquad |w|\, e^{-\gamma |w|} \le \frac{1}{e \gamma} .
\end{equation}
Indeed, $1 - e^{-w} = \int_0^w e^{-u}\, du$ with $|e^{-u}| \le 1$ on the segment of integration gives the upper bound. The functions $(1 - e^{-w}) / \min(1, |w|)$ and $1 + e^{-w}$ are continuous and zero-free on $\mathcal{W}$ (the zeros of $1 \mp e^{-w}$ lie on the imaginary axis), they have non-zero limits as $w \to 0$, and their moduli are at least $1 - e^{-\gamma |w|}$, which tends to $1$ as $|w| \to \infty$; this gives the lower bounds. The last inequality is elementary.

For $z \in \mathcal{D}(\eta_0, \phi)$ we have $|\arg(1 - ez)| < \pi - \phi$ and therefore $|\arg y| < \pi/2 - \phi/2$. By (T1), $\tau = y + \mathcal{O}(y^2)$; hence
\[
s := -\ln t = -\ln(1 - \tau) = y + \mathcal{O}(y^2)
\]
(principal branch) is well defined near $e^{-1}$, and there is $\delta_1 > 0$ such that for all $z \in \mathcal{D}(\eta_0, \phi)$ with $|z - e^{-1}| < \delta_1$:
\begin{equation}\label{eq:s}
s \in \mathcal{W}, \qquad |s| \le \tfrac{1}{2}, \qquad \tfrac{1}{2} |s| \le |\tau| \le |s| .
\end{equation}
For such $z$ we have $t^k = e^{-ks}$ with $ks \in \mathcal{W}$, so by \eqref{eq:sector}, for all $k \ge 1$,
\begin{equation}\label{eq:tk}
|t^k| \le e^{-\gamma k |s|}, \qquad \gamma' \min(1, k|s|) \le |1 - t^k| \le \min(2, k|s|), \qquad |1 + t^k| \ge \gamma' .
\end{equation}
We put
\[
\Lambda_k := \min\bigl( k+1,\ 1/|s| \bigr) \ge 1, \qquad \ell_k := 1 + \ln \Lambda_k .
\]
From \eqref{eq:tk} and the last inequality in \eqref{eq:sector},
\begin{equation}\label{eq:U}
\frac{|t^k|}{\Lambda_k} = |t^k| \max\Bigl( \frac{1}{k+1},\, |s| \Bigr) \le \frac{C_\gamma}{k+1} \qquad (k \ge 0), \qquad C_\gamma := 1 + \frac{2}{e\gamma} .
\end{equation}

Let $(d^*_h)$ be an admissible sequence, and suppose that at some point $z$ the values $d^*_0, \dots, d^*_{h-1}$ are non-zero and of modulus at most $1$. Then (b) can be written as $1/d^*_{j+1} = t^{-1} \bigl( 1/d^*_j + 1/2 + R(d^*_j) \bigr)$ for $0 \le j < h$. Multiplying by $t^{j+1}$ and summing over $j$ we obtain
\begin{equation}\label{eq:Ah}
A^*_h := \frac{t^h}{d^*_h} = \frac{1 - t^h}{2\tau} + b^*_h, \qquad b^*_h := \frac{1}{d^*_0} + \sum_{j=0}^{h-1} t^j R(d^*_j) .
\end{equation}
For the sequence $(d_h)$ itself we write $A_h$ and $b_h$ instead of $A^*_h$ and $b^*_h$.

\emph{Claim B. For every admissible sequence $(d^*_h)$ there exist $\delta \in (0, \delta_1]$ and constants $B, C_b \ge 1$, $c \in (0,1]$ such that for all $z$ in
\[
\mathcal{N} := \bigl\{ z : 0 < |z - e^{-1}| < \delta,\ |\arg(z - e^{-1})| > \phi \bigr\}
\]
and all $h \ge 0$:}
\begin{enumerate}
\item[(i)] $d^*_h \ne 0$ and $|d^*_h| \le \min\bigl(1, B/(h+1)\bigr)$;
\item[(ii)] $|b^*_h| \le C_b\, \ell_h$;
\item[(iii)] $c\, \Lambda_h \le |A^*_h| \le (2 + C_b)\, \Lambda_h$.
\end{enumerate}
Put $B := \max\bigl( 4,\ 8 C_\gamma / \gamma' \bigr)$ and $C_b := 2/x^*_0 + c_R B (3 + 2/\gamma)$, and choose an integer $N_0 \ge B$ such that $\frac{\gamma'}{8} \Lambda \ge C_b (1 + \ln \Lambda)$ for all $\Lambda \ge N_0$. As $z \to e^{-1}$ in $\mathcal{D}(\eta_0, \phi)$ we have $s \to 0$, $t \to 1$ and $d^*_h(z) \to x^*_h$ for every fixed $h$, where $x^*_h < \min(1, 2/(h+2))$ by \eqref{eq:xstar}. Hence we can choose $\delta \le \min(\delta_1, \eta_0 e^{-1})$ so small that on $\mathcal{N}$
\begin{equation}\label{eq:delta}
\frac{1}{|s|} \ge N_0, \qquad \frac{|s|}{1 - e^{-\gamma |s|}} \le \frac{2}{\gamma}, \qquad \Bigl| \frac{1}{d^*_0} \Bigr| \le \frac{2}{x^*_0}, \qquad |t|^{N_0} \ge \frac{1}{2},
\end{equation}
and such that $d^*_h \ne 0$ and $|d^*_h| \le \min\bigl(1, 4/(h+1)\bigr)$ on $\mathcal{N}$ for all $h < N_0$.
In particular (i) holds for $h < N_0$. Now let $h \ge 0$ and assume that (i) holds for all $j < h$. Then \eqref{eq:Ah} is valid, and by $|R(d^*_j)| \le c_R |d^*_j|$ and \eqref{eq:tk},
\begin{align*}
|b^*_h| \le \frac{2}{x^*_0} + c_R B \sum_{j=0}^{h-1} \frac{e^{-\gamma j |s|}}{j+1}
&\le \frac{2}{x^*_0} + c_R B \Bigl( \sum_{0 \le j < \Lambda_h} \frac{1}{j+1} + |s| \sum_{j \ge 1/|s|} e^{-\gamma j |s|} \Bigr)\\
&\le \frac{2}{x^*_0} + c_R B \Bigl( 2 + \ln \Lambda_h + \frac{2}{\gamma} \Bigr) \le C_b\, \ell_h ,
\end{align*}
which is (ii). If moreover $h \ge N_0$, then $\Lambda_h \ge N_0$, and by \eqref{eq:tk}, \eqref{eq:s} and the choice of $N_0$,
\[
|A^*_h| \ge \frac{|1 - t^h|}{2 |\tau|} - |b^*_h| \ge \frac{\gamma'}{2} \min\Bigl( h, \frac{1}{|s|} \Bigr) - C_b\, \ell_h \ge \frac{\gamma'}{4} \Lambda_h - C_b\, \ell_h \ge \frac{\gamma'}{8} \Lambda_h .
\]
Consequently $d^*_h = t^h / A^*_h \ne 0$ and, by \eqref{eq:U}, $|d^*_h| \le \frac{8}{\gamma'} \frac{|t^h|}{\Lambda_h} \le \frac{B}{h+1} \le 1$. By induction, (i) and (ii) hold for all $h \ge 0$. The upper bound in (iii) follows from \eqref{eq:Ah}, \eqref{eq:tk}, \eqref{eq:s} and $\ell_h \le \Lambda_h$:
\[
|A^*_h| \le \frac{\min(2, h|s|)}{|s|} + C_b\, \ell_h \le (2 + C_b) \Lambda_h .
\]
The lower bound in (iii) was shown above for $h \ge N_0$; for $h < N_0$ we have $|A^*_h| = |t|^h / |d^*_h| \ge 1/2 \ge \Lambda_h / (2 N_0)$. So (iii) holds with $c := \min(\gamma'/8,\ 1/(2N_0))$. This proves Claim~B.

From now on we apply Claim~B to the sequence $(d_h)$.

\emph{Claim C. After decreasing $\delta$ if necessary, there is a constant $C$ such that for all $z \in \mathcal{N}$ and $k \ge 0$:}
\begin{equation}\label{eq:claimC}
\frac{1}{C \Lambda_k} \le |\varepsilon_k| \le \frac{C}{\Lambda_k}, \qquad |d_k| \le \frac{C\, |t^k|}{\Lambda_k}, \qquad |\Delta_k| \le \frac{C\, |t^k|}{\Lambda_k^2} .
\end{equation}
By \eqref{eq:Ah}, $2 \tau A_k = 1 - t^k + 2\tau b_k$, and therefore
\[
\varepsilon_k = \tau + \frac{t^k}{A_k} = \frac{1 + t^k + 2 \tau b_k}{2 A_k} .
\]
By Claim~B(ii), $|2\tau b_k| \le 2 C_b |s| (1 + \ln(1/|s|))$, which is at most $\gamma'/2$ on $\mathcal{N}$ if $\delta$ is small enough. Together with $\gamma' \le |1 + t^k| \le 2$ this shows that the numerator has modulus between $\gamma'/2$ and $3$, and the bounds for $\varepsilon_k$ follow from Claim~B(iii). The bound for $d_k = t^k / A_k$ also follows from Claim~B(iii). For $k \ge 1$, \eqref{eq:rec} gives
\begin{equation}\label{eq:DeltaG}
\Delta_k = d_{k-1} - d_k = \tau\, d_{k-1} + t\, G(d_{k-1}),
\end{equation}
so $|\Delta_k| \le |s|\, |d_{k-1}| + |d_{k-1}|^2 \le C |t^{k-1}| / \Lambda_{k-1}^2$, because $|s| \le 1/\Lambda_{k-1}$; and $|t^{k-1}| \le 2 |t^k|$, $\Lambda_{k-1} \ge \Lambda_k / 2$. For $k = 0$ the bound is trivial. This proves Claim~C.

\paragraph{Proof of Lemma~\ref{lem:uniform}(1).} Let $\phi$, $\delta$ and $\mathcal{N}$ be as above, let $\eta \le \eta_0$, $\theta$ and $\mathcal{K}$ be given by Claim~A (applied to $(d_h)$) for this $\delta$, and put $\mathcal{D} := \mathcal{D}(\eta, \phi)$; then $\mathcal{D} \subseteq \mathcal{N} \cup \mathcal{K}$. On $\mathcal{D} \cap \mathcal{K}$ all assertions follow from Claim~A, since $|\Delta_k| \le |d_{k-1}| + |d_k| \le C \theta^{k}$ for $k \ge 1$. On $\mathcal{N}$ we have $\varepsilon_k \ne 0$ and $|t_k| \le |t| + |d_k| \le 2$ by Claims~B and~C. If $|k - k'| \le m$, then $\Lambda_{k'} \le (m+1) \Lambda_k$, so \eqref{eq:claimC} gives $|\varepsilon_k| \le C^2 (m+1) |\varepsilon_{k'}|$, and by \eqref{eq:claimC} and \eqref{eq:U}
\[
|\Delta_k| \le \frac{C\, |t^k|}{\Lambda_k} \cdot \frac{1}{\Lambda_k} \le \frac{C\, C_\gamma}{k+1} \cdot C\, |\varepsilon_k| \le \frac{C^4 C_\gamma (m+1)}{k+1}\, |\varepsilon_{k'}| .
\]

\paragraph{Step 3: the sum.} We now prove Lemma~\ref{lem:uniform}(2). Fix $a, b, p, q, r \in \{0, \dots, m\}$ and let $T_k(z)$ denote the $k$-th term of the series. In the sequel $C$ denotes constants that may depend on $m$ and $\phi$ and may change from line to line.

\emph{Convergence.} On $\mathcal{D} \cap \mathcal{K}$ we have $|T_k| \le C \theta^{2k}$ by Claim~A. On $\mathcal{N}$, \eqref{eq:claimC} and \eqref{eq:tk} give $|T_k| \le C |t^k|^2 / \Lambda_k \le C e^{-2\gamma k |s|}$, and $|s|$ is bounded away from $0$ on compact subsets of $\mathcal{N}$. Hence the series converges locally uniformly in $\mathcal{D}$.

\emph{The model.} For $z \in \mathcal{N}$ and $k \ge 1$ define
\[
\bar d_k := \frac{2 \tau\, t^k}{1 - t^k}, \qquad \bar\varepsilon_k := \tau + \bar d_k = \tau\, \frac{1 + t^k}{1 - t^k}, \qquad \bar\Delta_k := \tau\, \bar d_k + \tfrac{1}{2} \bar d_k^{\,2} = \frac{2 \tau^2\, t^k}{(1 - t^k)^2} .
\]
By \eqref{eq:tk} and \eqref{eq:s}, $\frac{1}{4 \Lambda_k} \le \frac{|\tau|}{|1 - t^k|} \le \frac{2}{\gamma' \Lambda_k}$, hence
\begin{equation}\label{eq:modelsize}
|\bar d_k| \le \frac{C\, |t^k|}{\Lambda_k}, \qquad \frac{1}{C \Lambda_k} \le |\bar\varepsilon_k| \le \frac{C}{\Lambda_k}, \qquad |\bar\Delta_k| \le \frac{C\, |t^k|}{\Lambda_k^2} .
\end{equation}

\emph{Claim D. For all $z \in \mathcal{N}$, $k \ge 1$ and $0 \le l \le m$:}
\[
|\varepsilon_{k+l} - \bar\varepsilon_k| \le \frac{C\, |t^k|\, \ell_k}{\Lambda_k^2}, \qquad |\Delta_{k+l} - \bar\Delta_k| \le \frac{C\, |t^k|\, \ell_k}{\Lambda_k^3} .
\]
For $j \ge 1$, \eqref{eq:Ah} gives
\[
d_j - \bar d_j = t^j \Bigl( \frac{1}{A_j} - \frac{2\tau}{1 - t^j} \Bigr) = - \frac{t^j\, b_j}{A_j \cdot (1 - t^j)/(2\tau)},
\]
so $|d_j - \bar d_j| \le C |t^j| \ell_j / \Lambda_j^2$ by Claim~B. Next, let $j = k + l'$ with $-1 \le l' \le m$ and $j \ge 1$. Then
\[
\bar d_j - \bar d_k = \frac{2 \tau\, t^k\, (t^{l'} - 1)}{(1 - t^j)(1 - t^k)}, \qquad |t^{l'} - 1| \le C |s|,
\]
so $|\bar d_j - \bar d_k| \le C |t^k| / \Lambda_k^2$. Since $|t^j| \le 2 |t^k|$, $\Lambda_j \ge \Lambda_k / 2$ and $\ell_j \le C \ell_k$, we conclude that
\begin{equation}\label{eq:dshift}
|d_{k+l'} - \bar d_k| \le \frac{C\, |t^k|\, \ell_k}{\Lambda_k^2} \qquad (k \ge 1,\ -1 \le l' \le m);
\end{equation}
the only case not covered so far is $k = 1$, $l' = -1$, in which the left-hand side is bounded and the right-hand side is bounded away from $0$. The first assertion of Claim~D follows, as $\varepsilon_{k+l} - \bar\varepsilon_k = d_{k+l} - \bar d_k$. For the second assertion put $j := k + l - 1$. By \eqref{eq:DeltaG},
\[
\Delta_{k+l} - \bar\Delta_k = \tau (d_j - \bar d_k) + t \bigl( G(d_j) - \tfrac{1}{2} d_j^2 \bigr) - \tfrac{1}{2} \tau\, d_j^2 + \tfrac{1}{2} (d_j - \bar d_k)(d_j + \bar d_k) .
\]
Using $|\tau| \le 1/\Lambda_k$, $|G(x) - x^2/2| \le |x|^3$, \eqref{eq:dshift}, and the bounds $|d_j|, |\bar d_k| \le C |t^k| / \Lambda_k$ from \eqref{eq:claimC} and \eqref{eq:modelsize}, each of the four terms is at most $C |t^k| \ell_k / \Lambda_k^3$. This proves Claim~D.

\emph{Comparison with the model.} For $k \ge 1$ let
\[
P_k := \frac{\bar\Delta_k^{\,2}}{\bar\varepsilon_k^{\,3}} = \frac{4 \tau\, t^{2k}}{(1 - t^k)(1 + t^k)^3} .
\]
Replacing the five factors of $T_k$ by those of $P_k$ one at a time and using \eqref{eq:claimC}, \eqref{eq:modelsize}, Claim~D and then \eqref{eq:U}, we get for $z \in \mathcal{N}$ and $k \ge 1$
\begin{equation}\label{eq:dom}
|T_k - P_k| \le C\, \Bigl( \frac{|t^k|}{\Lambda_k} \Bigr)^{2} \ell_k \le \frac{C\, (1 + \ln(k+1))}{(k+1)^2} .
\end{equation}
The right-hand side is summable and does not depend on $z$. Moreover, as $z \to e^{-1}$ in $\mathcal{D}$, each $T_k(z)$ tends to $T_k(e^{-1})$, which is finite because $\varepsilon_h(e^{-1}) = x_h > 0$, and $P_k(z) \to \frac{1}{2k}$ because $(1 - t^k)/\tau \to k$. By dominated convergence,
\begin{equation}\label{eq:kappa1}
\sum_{k \ge 1} \bigl( T_k(z) - P_k(z) \bigr) = \kappa_1 + o(1), \qquad \kappa_1 := \sum_{k \ge 1} \Bigl( T_k(e^{-1}) - \frac{1}{2k} \Bigr).
\end{equation}

\emph{The model sum.} Let
\[
\varphi(w) := \frac{4 e^{-2w}}{(1 - e^{-w})(1 + e^{-w})^3}, \qquad \chi(w) := \varphi(w) - \frac{e^{-w}}{2w},
\]
so that $P_k = \tau\, \varphi(ks)$. The function $\varphi$ has a simple pole with residue $\frac{1}{2}$ at $w = 0$, so $\chi$ is analytic in a neighbourhood of $\mathcal{W} \cup \{0\}$. By \eqref{eq:sector}, applied in a slightly larger sector, and Cauchy's estimates, there are constants $C$ and $\gamma_1 > 0$ with $|\chi(w)| + |\chi'(w)| \le C e^{-\gamma_1 |w|}$ for $w \in \mathcal{W}$. Since $\sum_{k \ge 1} e^{-ks}/k = -\ln(1 - e^{-s}) = \ln(1/\tau)$,
\[
\sum_{k \ge 1} P_k = \frac{\tau}{s} \sum_{k \ge 1} s\, \chi(ks) + \frac{\tau}{2s} \ln \frac{1}{\tau} .
\]
The first sum is a Riemann sum:
\[
\Bigl| \sum_{k \ge 1} s\, \chi(ks) - \int_0^\infty s\, \chi(xs)\, dx \Bigr| \le |s| \sum_{k \ge 1} \int_{k-1}^{k} |\chi(ks) - \chi(xs)|\, dx \le |s|^2 \int_0^\infty |\chi'(xs)|\, dx \le C |s|,
\]
and $\int_0^\infty s\, \chi(xs)\, dx = \int_0^\infty \chi(x)\, dx =: I$ by Cauchy's theorem, the integrand being analytic and exponentially small in $\mathcal{W}$. Since $\tau / s = 1 + \mathcal{O}(|s|)$, we obtain
\[
\sum_{k \ge 1} P_k = \frac{1}{2} \ln \frac{1}{\tau} + I + \mathcal{O}\bigl( |s| \ln (1/|s|) \bigr) .
\]
Finally $\tau = y (1 + \mathcal{O}(y))$ and $y^2 = 2(1 - ez)$, so $\ln(1/\tau) = \frac{1}{2} L(z) - \frac{1}{2} \ln 2 + \mathcal{O}(|y|)$. Together with \eqref{eq:kappa1} and $T_0(z) \to T_0(e^{-1})$ this gives
\[
\sum_{k \ge 0} T_k(z) = \frac{1}{4} L(z) + \kappa + o(1), \qquad \kappa := T_0(e^{-1}) + \kappa_1 + I - \frac{1}{4} \ln 2,
\]
as $z \to e^{-1}$ in $\mathcal{D}$. This completes the proof of Lemma~\ref{lem:uniform}. \qed

\section{Proof of Proposition~\ref{prop:Mhat}}
\label{app:crown}

We use the notation of Section~\ref{sec:crown} and of Appendix~\ref{app:uniform}; in particular, $c \ge 0$, $\alpha > 1$ and $\lambda$ are fixed as in Section~\ref{sec:crown}.

\paragraph{Step 1: the tilted sequence.} Put $\tilde d_k := t - \tilde t_k$ for $k \ge 0$, so that
\[
\hat Y_k = d_k - \tilde d_k .
\]
Since $\tilde t_{k+1} = z e^{\tilde t_k}$, the sequence $(\tilde d_k)$ satisfies the recurrence \eqref{eq:rec}, and $\tilde d_0 = t - z e^{S} = t \bigl( 1 - e^{S - t} \bigr)$ because $z = t e^{-t}$. We check that $(\tilde d_k)$ is admissible in the sense of Appendix~\ref{app:uniform}. Conditions (a) and (b) hold (we may assume that $1 + \eta_0 < \lambda$). As $z \to e^{-1}$ we have $\tilde d_0 \to \tilde x_0 := 1 - e^{S_0 - 1} \in (0,1)$ by \eqref{eq:S0}, which is (d); consequently $\tilde x_k := \tilde d_k(e^{-1}) > 0$ for all $k$. Finally, by \eqref{eq:Yk} and \eqref{eq:inside}, for $|z| \le e^{-1}$
\[
|\tilde d_k(z)| \le |d_k(z)| + |\hat Y_k(z)| \le x_k + \hat Y_k(e^{-1}) = 2 x_k - \tilde x_k < 2 x_k ,
\]
which tends to $0$ by \eqref{eq:xh}; this is (c).

Fix $\phi \in (0, \pi/2)$. We use the notation and the results of Appendix~\ref{app:uniform}, and we choose $\delta$ and $\eta$ so small that, with $\mathcal{N}$, $\mathcal{K}$ and $\mathcal{D} := \mathcal{D}(\eta, \phi) \subseteq \mathcal{N} \cup \mathcal{K}$ as there, Claims~A and~B hold for both sequences $(d_k)$ and $(\tilde d_k)$, Claims~C and~D hold with $m := c+2$, and Proposition~\ref{prop:gap_sing} holds in $\mathcal{D}$ for $i = 1$. We write $\tilde A_k$, $\tilde b_k$ for the quantities \eqref{eq:Ah} belonging to $(\tilde d_k)$, and
\[
\beta_k := \tilde b_k - b_k = \tilde A_k - A_k .
\]
In the sequel $C$ denotes constants that do not depend on $z$ and $k$. On $\mathcal{N}$ we have, by \eqref{eq:Ah},
\begin{equation}\label{eq:Yhat}
\hat Y_k = d_k - \tilde d_k = t^k \Bigl( \frac{1}{A_k} - \frac{1}{\tilde A_k} \Bigr) = \frac{t^k\, \beta_k}{A_k\, \tilde A_k} .
\end{equation}
The function $R$ is Lipschitz on the unit disc, and $|d_j|, |\tilde d_j| \le 1$ by Claim~B(i). Hence, by \eqref{eq:Yhat}, Claim~B(iii) and \eqref{eq:U},
\begin{equation}\label{eq:Rdiff}
\bigl| t^j \bigl( R(\tilde d_j) - R(d_j) \bigr) \bigr| \le C\, |t^j|\, |\hat Y_j| \le C \Bigl( \frac{|t^j|}{\Lambda_j} \Bigr)^{2} |\beta_j| \le \frac{C_3\, |\beta_j|}{(j+1)^2}
\end{equation}
for a constant $C_3$. Since $\beta_k = \beta_0 + \sum_{j<k} t^j \bigl( R(\tilde d_j) - R(d_j) \bigr)$ and $|\beta_0| = |1/\tilde d_0 - 1/d_0| \le 2/\tilde x_0 + 2/x_0$ on $\mathcal{N}$ by \eqref{eq:delta}, induction on $k$ gives
\[
|\beta_k| \le |\beta_0| \prod_{j<k} \Bigl( 1 + \frac{C_3}{(j+1)^2} \Bigr) \le C \qquad (z \in \mathcal{N},\ k \ge 0).
\]
Inserting this into \eqref{eq:Rdiff} we see that $\beta_k$ converges uniformly on $\mathcal{N}$ to an analytic function $\beta_\infty$, and
\begin{equation}\label{eq:betatail}
|\beta_\infty(z) - \beta_k(z)| \le \frac{C}{k+1}, \qquad |\beta_\infty(z)| \le C \qquad (z \in \mathcal{N},\ k \ge 0).
\end{equation}

\paragraph{Step 2: $\hat M$ near $e^{-1}$.} We first note that the series $\hat M$ converges locally uniformly in $\mathcal{D}$, and that $\hat M$ is bounded on $\mathcal{D} \cap \mathcal{K}$: there $|\hat Y_k| \le |d_k| + |\tilde d_k| \le C \theta^k$ and $|\Psi_k| \le C$ by Claim~A. On $\mathcal{N}$ we have, for $k \ge 1$,
\begin{equation}\label{eq:Psibound}
|\Psi_k| \le C\, |t^k|\, \Lambda_k ,
\end{equation}
which follows from \eqref{eq:psidiff}, \eqref{eq:claimC}, $|t_h| \le 2$ and $t_{k+c} - t_k = \sum_{l=1}^{c} \Delta_{k+l}$; together with \eqref{eq:Yhat} and Claim~B(iii) this gives $|\hat Y_k \Psi_k| \le C |t^k|^2 / \Lambda_k \le C e^{-2\gamma k |s|}$, and local uniform convergence follows as in Step~3 of Appendix~\ref{app:uniform}. In particular, $\hat M$ is analytic in $\mathcal{D}$, and by Weierstrass' theorem it is the analytic continuation of the power series $\hat M(z)$.

We now show that
\begin{equation}\label{eq:Mhat}
\hat M(z) = \beta_\infty(z)\, L(z) + \mathcal{O}(1) \qquad (z \in \mathcal{N}).
\end{equation}
Let $k \ge 1$ and recall $\bar\Delta_k = 2 \tau^2 t^k / (1 - t^k)^2$. As in the proof of Claim~D, Claim~B gives
\[
\Bigl| \frac{1}{A_k} - \frac{2\tau}{1 - t^k} \Bigr| \le \frac{C\, \ell_k}{\Lambda_k^2}, \qquad \Bigl| \frac{1}{\tilde A_k} - \frac{2\tau}{1 - t^k} \Bigr| \le \frac{C\, \ell_k}{\Lambda_k^2},
\]
and all three quantities $1/A_k$, $1/\tilde A_k$, $2\tau/(1 - t^k)$ are bounded by $C / \Lambda_k$. Hence $\bigl| t^k / (A_k \tilde A_k) - 2 \bar\Delta_k \bigr| \le C |t^k| \ell_k / \Lambda_k^3$, and by \eqref{eq:Yhat}, \eqref{eq:betatail} and \eqref{eq:modelsize},
\[
\bigl| \hat Y_k - 2 \beta_\infty \bar\Delta_k \bigr| \le |\beta_k|\, \Bigl| \frac{t^k}{A_k \tilde A_k} - 2 \bar\Delta_k \Bigr| + 2\, |\beta_k - \beta_\infty|\, |\bar\Delta_k| \le C\, |t^k| \Bigl( \frac{\ell_k}{\Lambda_k^3} + \frac{1}{(k+1) \Lambda_k^2} \Bigr).
\]
By Claim~D, $|\Delta_{k+1} - \bar\Delta_k| \le C |t^k| \ell_k / \Lambda_k^3$. Combining these two estimates with \eqref{eq:Psibound} and \eqref{eq:U} we obtain
\[
\bigl| \bigl( \hat Y_k - 2 \beta_\infty \Delta_{k+1} \bigr) \Psi_k \bigr| \le C \Bigl( \frac{|t^k|}{\Lambda_k} \Bigr)^{2} \ell_k + \frac{C\, |t^k|^2}{(k+1) \Lambda_k} \le \frac{C\, (1 + \ln(k+1))}{(k+1)^2} \qquad (k \ge 1).
\]
The right-hand side is summable, and the terms with $k = 0$ are bounded on $\mathcal{N}$. Since $\sum_{k \ge 0} \Delta_{k+1} \Psi_k = E^{(c)}_1(z)$ by \eqref{eq:Ei}, it follows that
\[
\hat M(z) = 2 \beta_\infty(z)\, E^{(c)}_1(z) + \mathcal{O}(1) \qquad (z \in \mathcal{N}),
\]
and \eqref{eq:Mhat} follows from Proposition~\ref{prop:gap_sing} and the boundedness of $\beta_\infty$.

\paragraph{Step 3: $\beta_\infty$ near $e^{-1}$.} We show that there is a constant $\beta^*$ such that, after decreasing $\delta$ if necessary,
\begin{equation}\label{eq:betarate}
|\beta_\infty(z) - \beta^*| \le \frac{C}{\ln(1/|y|)} \qquad (z \in \mathcal{N}),
\end{equation}
where $y = \sqrt{2(1-ez)}$ as before. Recall that the tree function is an analytic function of $y$: there is a function $\mathfrak{t}$, analytic in a disc $|y| < y_0$, with $\mathfrak{t}(0) = 1$ and $t(z) = \mathfrak{t}(y)$ for $z \in \mathcal{D}(\eta_0, \phi)$ close to $e^{-1}$ \cite{flajolet2009}. (Indeed, $z = t e^{-t}$ gives $y^2 = 2(1 - t e^{1-t}) = (1-t)^2 (1 + \mathcal{O}(1-t))$, so $y$ is an analytic function of $1 - t$ with derivative $1$ at the origin, which can be inverted.) With $z = (1 - y^2/2) e^{-1}$, the functions $d_0 = \mathfrak{t} - z$ and $\tilde d_0 = \mathfrak{t} (1 - e^{S(z) - \mathfrak{t}})$, and hence, by the recurrence, all $d_j$ and $\tilde d_j$, become analytic functions of $y$ in the disc $|y| < y_0$, with the values $x_j$ and $\tilde x_j$ at $y = 0$. We use the same letters for these functions of $y$; note that now $y$ ranges over a full disc and not only over the sector corresponding to $\mathcal{N}$.

Choose $y_1 \in (0, y_0)$ and $C_0 \ge 8$ such that for $|y| \le y_1$
\[
|\mathfrak{t}(y) - 1| \le 2 |y|, \qquad |d_0(y) - x_0| \le C_0 |y|, \qquad |\tilde d_0(y) - \tilde x_0| \le C_0 |y| .
\]
We claim that for all $j \ge 0$ and all $y$ with $|y| \le y_1$ and $C_0\, 3^j |y| \le \frac{1}{2}$
\begin{equation}\label{eq:lip}
|d_j(y) - x_j| \le C_0\, 3^j |y| \qquad \text{and} \qquad |\tilde d_j(y) - \tilde x_j| \le C_0\, 3^j |y| .
\end{equation}
This is true for $j = 0$. Assume it for $j$ and let $C_0 3^{j+1} |y| \le \frac{1}{2}$. Then $w := d_j(y) - x_j$ satisfies $|w| \le \frac{1}{6}$, and
\[
d_{j+1}(y) - x_{j+1} = (\mathfrak{t} - 1) \bigl( 1 - e^{-d_j(y)} \bigr) + e^{-x_j} \bigl( 1 - e^{-w} \bigr) .
\]
Here $|1 - e^{-d_j(y)}| \le 1 + e^{|w|} \le 3$ and $|1 - e^{-w}| \le |w| e^{|w|} \le 2|w|$, so $|d_{j+1}(y) - x_{j+1}| \le 6 |y| + 2 C_0 3^j |y| \le C_0 3^{j+1} |y|$. The same argument applies to $\tilde d_j$, and \eqref{eq:lip} follows.

Since $1/x_{h+1} = 1/x_h + 1/2 + R(x_h) \le 1/x_h + 1/2 + c_R$, and similarly for $\tilde x_h$, there is a constant $C_1$ with
\[
x_K \ge \frac{1}{C_1 (K+1)}, \qquad \tilde x_K \ge \frac{1}{C_1 (K+1)} \qquad (K \ge 0).
\]
Choose $C_2$ so large that $\varrho_K := 4^{-K} / C_2$ satisfies $\varrho_K \le y_1$, $(1 + 2 \varrho_K)^K \le 2$ and $C_0 3^K \varrho_K \le \frac{1}{2 C_1 (K+1)}$ for all $K \ge 0$. By \eqref{eq:lip}, for $|y| \le \varrho_K$ we have $|d_K(y)| \ge x_K / 2$ and $|\tilde d_K(y)| \ge \tilde x_K / 2$. Hence the function
\[
\beta_K(y) := \mathfrak{t}(y)^K \Bigl( \frac{1}{\tilde d_K(y)} - \frac{1}{d_K(y)} \Bigr)
\]
is analytic for $|y| \le \varrho_K$ and bounded there by $8 C_1 (K+1)$. By the definition of $A_K$ and $\tilde A_K$ in \eqref{eq:Ah}, it coincides with the function $\beta_K = \tilde A_K - A_K$ of Step~1 at the points $y$ that correspond to $z \in \mathcal{N}$. By the Schwarz lemma,
\begin{equation}\label{eq:schwarz}
|\beta_K(y) - \beta_K(0)| \le 16\, C_1 (K+1)\, \frac{|y|}{\varrho_K} = 16\, C_1 C_2\, (K+1)\, 4^K |y| \qquad (|y| \le \varrho_K).
\end{equation}
Letting $z \to e^{-1}$ in \eqref{eq:betatail} we get $|\beta_{K'}(0) - \beta_K(0)| \le 2C/(K+1)$ for $K' \ge K$; hence $\beta^* := \lim_{K \to \infty} \beta_K(0)$ exists and $|\beta^* - \beta_K(0)| \le 2C/(K+1)$.

Now let $z \in \mathcal{N}$ and put $K := \lfloor \ln(1/|y|) / \ln 8 \rfloor$. Then $|y| \le 8^{-K}$, and if $\delta$ is small enough, $K$ is so large that $8^{-K} \le \varrho_K$. By \eqref{eq:schwarz}, $|\beta_K(y) - \beta_K(0)| \le 16\, C_1 C_2 (K+1) 2^{-K} \le C/(K+1)$. Together with \eqref{eq:betatail} and the bound for $\beta^* - \beta_K(0)$ this yields
\[
|\beta_\infty(z) - \beta^*| \le \frac{C}{K+1} \le \frac{C \ln 8}{\ln(1/|y|)},
\]
which is \eqref{eq:betarate}.

\paragraph{Step 4: conclusion.} Since $|L(z)| = |\ln(2/y^2)| \le 2 \ln(1/|y|) + \ln 2 + \pi$, the estimates \eqref{eq:Mhat} and \eqref{eq:betarate} give
\[
\hat M(z) = \beta^* L(z) + \mathcal{O}(1) \qquad (z \in \mathcal{N}).
\]
As $\hat M$ and $L$ are bounded on $\mathcal{D} \cap \mathcal{K}$, the function $\hat M - \beta^* L$ is bounded in the $\Delta$-domain $\mathcal{D}$. This proves Proposition~\ref{prop:Mhat}. \qed

\bibliographystyle{plain}
\bibliography{htree}

\end{document}